\documentclass[review]{elsarticle}

\usepackage{lineno,hyperref}
\usepackage{amsmath}
\usepackage{amssymb}
\usepackage{bm}
\usepackage{subcaption}
\usepackage{xcolor}
\modulolinenumbers[5]
\usepackage{algorithm}
\usepackage{algorithmic}
\usepackage{multirow}
\usepackage{diagbox}
\journal{arXiv}

\newcommand{\ksi}{\xi}
\newcommand{\vksi}{\boldsymbol{\ksi}}

\def\bu{\boldsymbol{u}}

\def\pvt{ \boldsymbol{\mathcal{P}}}

\def\pp{\partial}

\def\bu{\boldsymbol{u}}

\def\R{\boldsymbol{R}}

\newcommand{\herm}{\Psi}

\def\bl{\boldsymbol{\lambda}}

\begin{document}
\nolinenumbers
\begin{frontmatter}

\title{Sensitivity-enhanced generalized polynomial chaos for efficient uncertainty quantification}

\author{Kyriakos D. Kantarakias}
\author{George Papadakis}
\address{Department of Aeronautics, Imperial College London, SW7 2AZ, U.K.}

\begin{abstract}
We present an enriched formulation of the Least Squares (LSQ) regression method for Uncertainty Quantification (UQ) using generalised polynomial chaos (gPC). More specifically, we enrich the linear system with additional equations for the gradient (or sensitivity) of the Quantity of Interest with respect to the stochastic variables. This sensitivity is computed very efficiently for all variables by solving an adjoint system of equations at each sampling point of the stochastic space. The associated computational cost is similar to one solution of the direct problem. For the selection of the sampling points, we apply a greedy algorithm which is based on the pivoted QR decomposition of the measurement matrix. We call the new approach sensitivity-enhanced generalised polynomial chaos, or se-gPC. We apply the method to several test cases to test accuracy and convergence with increasing chaos order, including an aerodynamic case with $40$ stochastic parameters. The method is found to produce accurate estimations of the statistical moments using the minimum number of sampling points. The computational cost scales as $\sim m^{p-1}$, instead of $\sim m^p$ of the standard LSQ formulation, where $m$ is the number of stochastic variables and $p$ the chaos order. The solution of the adjoint system of equations is implemented in many computational mechanics packages, thus the infrastructure exists for the application of the method to a wide variety of engineering problems. 

\end{abstract}

\begin{keyword}
\texttt{elsarticle.cls}\sep \LaTeX\sep Elsevier \sep template
\MSC[2010] 00-01\sep  99-00
\end{keyword}

\end{frontmatter}

%\linenumbers
%%%%%%%%%%%%%%%%%%%%%%%%%%%%%%%%%%%%%%%%%%%%%%%%%%%%%%%%%%%%%%%%%%%%%%%
%%%%%%%%%%%%%%%%%%%%%%%%%%%%%%%%%%%%%%%%%%%%%%%%%%%%%%%%%%
\section{Introduction}
%%%%%%%%%%%%%%%%%%%%%%%%%%%%%%%%%%%%%%%%%%%%%%%%%%%%%%%%%%
\label{sec:Intro}
In all practical applications, the performance of an engineering system is affected by stochastic uncertainties in the system parameters, boundary conditions etc.  Significant research has been directed towards quantifying the effect of such uncertainties on the performance of engineering systems. The performance is evaluated by a metric, called Quantity of Interest (QoI), that depends on the system under investigation. For example, in the area of aerodynamics it can be the lift or drag coefficient of an airfoil, and the uncertainty may arise due to stochastic variations in the angle of attack, approaching velocity etc. This area of research is known as Uncertainty Quantification (UQ). 

The simplest UQ approaches are random sampling techniques, such as the Monte--Carlo (MC) \cite{MCbook07,MOROKOFF1995218}. Such  approaches however require a large number of samples for the statistics of the QoI to converge and are therefore computationally expensive. On the other hand, methods based on Polynomial Chaos Expansion (PCE) have much reduced computational cost.  In these approaches, the uncertain QoI is expanded in terms of multivariate polynomial bases that are orthogonal with respect to the joint probability density function (PDF) of the uncertain parameters. The method was first introduced in \cite{homoChaos} for Gaussian stochastic variables and later extended for other distributions in \cite{KarniadakisOriginal}. These methods have been very successful in predicting the statistics of the QoI in a wide range of applications, such as applied mechanics, fluid dynamics, space, medicine. etc \cite{GhanemSpanos, CHASSAING2012394,SCHIAVAZZI2017196,CHATZIMANOLAKIS2019207,PREUQKantar,JONES20131860}. 

In the standard Galerkin-projection formulation of these methods, the cost of computing the spectral coefficients of the QoI expansion scales exponentially with the dimension of the stochastic space. For the  non-intrusive projection methods, this is due to the large number of sampling points when extending one-dimensional quadrature integration schemes to multi-dimensions using tensor products, for details see the book of \cite{KnioBook}. At each sampling point, the evaluation of the model describing the system is by far the most time consuming part of the method. In sparse grid approaches such as Smolyak grids,  \cite{SmolyakPaper,GANAPATHYSUBRAMANIAN2007652,CONSTANTINE20121,sparse2}, the spectral coefficients are computed from evaluations at a significantly smaller number of sample points. In \cite{BLATMAN20112345} an adaptive method to build a sparse polynomial chaos (PC) expansion basis using  least-angle regression (LARS) is introduced that requires a reduced number of evaluations. 

Another approach to compute the spectral coefficients is based on least-squares (LSQ) regression, that minimises the expectation of the square of the difference between the QoI evaluation and the PCE prediction at different sampling points,  \cite{HAMPTON2015363,BLATMAN20112345,KnioBook}. A review of the techniques used to select the evaluation sampling points is provided in \cite{HADIGOL2018382}. Several optimality criteria are discussed, such as the alphabetic D-, A-, E-, K- and I-criteria, and a hybrid sampling method, called alphabetic-coherence-optimal, is proposed. 

For the least square regression to provide accurate results, the model has to be evaluated at a number of sampling points, say $q$, which scales with the number of unknown PCE coefficients. The oversampling ratio is usually between 2-3 \cite{BLATMAN20112345} or larger \cite{HADIGOL2018382}. For $m$ stochastic variables and polynomial chaos order $p$, the number of PCE coefficients is $\frac{(p+m)!}{m!p!}$, see also equation \eqref{theory03} below, which behaves asymptotically as $\sim \frac{1}{p!}m^p$ when $m\to \infty$, \cite{BLATMAN20112345}. Thus, the number of evaluations grows rapidly with increasing $p$ and $m$, and the computational cost quickly becomes intractable. Similarly, for Smolyak quadrature the computational cost scales as $\sim \frac{2^p}{p!}m^p$, \cite{BLATMAN20112345}. For the smallest chaos order $p=1$, the cost grows linearly with $m$ in either approach.

In the present paper, we aim to reduce the computational cost of the LSQ method. In the standard formulation, each model evaluation provides one equation for the linear system of the PCE coefficients.  We reformulate the method by enriching the system with additional equations that can be obtained very efficiently. More specifically, if the underlying system is described by a set of linear or non-linear Partial Differential Equations (PDE), the solution of the adjoint system allows for a very efficient evaluation of the sensitivity of the QoI with respect to all $m$ stochastic variables at each sampling point, \cite{LagrangeMultipliers, Luchini2014AdjointAnalysis, adjoint1,adjointjameson}. The cost of the solution of the (linear) adjoint system is comparable to that of the governing PDE. Therefore at each sampling point we can obtain $1+m$ equations at the cost of two evaluations, one for the direct and one for the adjoint system. Thus for $q$ sampling points, we have a linear system of $q(1+m)$ equations, and if we make this number of equations commensurate to the number of unknown PCE coefficients, we get  $q(1+m) \sim \frac{(p+m)!}{m!p!}$, which for large $m$ becomes $q(1+m) \sim m^p$ or $q \sim m^{p-1}$. Thus the computational cost now scales as $\sim m^{p-1}$, as opposed to $\sim m^p$ of the standard formulation. For example, for $p=1$ the cost becomes independent of $m$. This drastic reduction of the computational cost, by a factor $m$, is the central contribution of this paper. 

The solution of the adjoint equations is implemented in many computational mechanics packages (commercial or open source). This sensitivity information is mainly used for optimisation purposes (using gradient descent for example), but it can be also employed for UQ as described in the rest of the paper. In order to minimise the number of sampling points, and make the oversampling ratio close to 1, we employ the pivoted QR decomposition of the measurement matrix \cite{SOMMARIVA20091324,  BruntonPaper,DIAZ2018640}. We call the new approach sensitivity-enhanced gPC, or se-gPC.

The paper is organized as follows. In section \ref{sec:UQIntro}, an overview of UQ using gPC with LSQ is outlined. Section \ref{segPC basics} presents the formulation of sensitivity-enhanced gPC, followed by section \ref{sec:adjoint basics} that reviews the adjoint method for efficient computation of sensitivities. Section \ref{sec:Fekete Sampling} discusses the sampling of the input stochastic space. Results for several test cases are presented in section \ref{sec:High Dimensional Systems} and we conclude in section \ref{sec:conclusions chapter}. 

%%%%%%%%%%%%%%%%%%%%%%%%%%%%%%%%%%%%%%%%%%%%%%%%%%%%%%%%%%
\section{Uncertainty Quantification with  gPC}\label{sec:UQIntro}
%%%%%%%%%%%%%%%%%%%%%%%%%%%%%%%%%%%%%%%%%%%%%%%%%%%%%%%%%%
Let $\mathbb{P} = (\Omega, \Sigma,d\mathcal{P})$ denote a complete probability space, where $\Omega$ is the sample set of random events, and $d\mathcal{P}$ is a probability measure defined on the $\sigma$-algebra $\Sigma$. An uncertain input is modelled as a stochastic vector $\vksi$ of $m$ independent stochastic variables $\ksi_i$, i.e.\ $\vksi = [ \ksi_1, \dots, \ksi_m] : \Omega \to \mathbb{R}^m$ with $m \in \mathbb{N}$. We use subscripts to denote random variables, and superscripts for realisations, for example the $i-th$ realisation of vector $\vksi$ is $\vksi^{(i)} = [ \ksi_1^{(i)}, \dots, \ksi_m^{(i)}]$. The stochastic variable $\ksi_i$ follows a probability density function (PDF) $w_i(\ksi_i)$, defined in the domain $\mathcal{E}_i$. Since the variables are independent, the $m$-dimensional vector $\vksi$ follows the PDF $W = \prod_{i=1}^m w_i({\ksi_i})$, defined in the domain $\mathcal{E} = \prod_{i=1}^m \mathcal{E}_i $.
A polynomial basis with members $\Psi_j(\vksi)$ that is orthogonal with respect to $W$ is defined as,
\begin{equation}
\langle \Psi_j,\Psi_k \rangle = \int_{\mathcal{E}} \Psi_j \Psi_k W d\vksi = \delta_{jk } \langle \Psi_j,\Psi_j \rangle. 
\label{theory01}
\end{equation}
The polynomials are normalized so that $\left <\Psi_j,\Psi_j \right > = 1$. Note that unless specified otherwise, repeated indices don't imply summation. When $m>1$, the  multidimensional basis functions $\Psi_j$ are defined as tensor products of the univariate polynomials $\psi$, for details see \cite{KnioBook}. 

A Quantity of Interest (QoI) ${\mathcal{M}: \vksi \to \mathcal{M} \in \mathbb{R}} $ with finite variance can be written in terms of  $\herm_i$ as,
\begin{equation}
\mathcal{M} (\vksi) =\sum_{i=0}^{\infty} c^{i} \herm_i(\vksi)
\label{eq:pce02}
\end{equation}

In practise, the above expansion is truncated by limiting the maximum order of the polynomials $\herm_i(\vksi)$ to a value $p$, known as chaos order. In this case, we get 
\begin{equation}
\mathcal{M} (\vksi) =\sum_{i=0}^{P} c^{i} \herm_i(\vksi) + \epsilon(\vksi), 
\label{eq:pce02_truncated}
\end{equation}
where $\epsilon(\vksi)$ is the truncation error, and the number of retained $P+1$ terms is given by
\begin{equation}
P+1 = \frac{(p+m)!}{p! m!}.
\label{theory03}
\end{equation}
The statistical moments of $\mathcal{M}$ can be computed algebraically from the spectral coefficients $c^i$, 
\begin{equation}
\begin{aligned}
&\mbox{Mean value:} \quad E[\mathcal{M}] = c^0 \quad \quad \mbox{Variance:} \quad \sigma^2[\mathcal{M}] = \sum_{i=1}^P c^i\\
&\mbox{Skewness:} \quad S[\mathcal{M}] = \frac{E[\mathcal{M}^3] - 3E[\mathcal{M}] \sigma^2[\mathcal{M}] - E^3[\mathcal{M}]}{ \sigma^3[\mathcal{M}]} \\ 
& \mbox{Kurtosis:} \quad K[\mathcal{M}] = \frac{E[\mathcal{M}^4] - 4 E[\mathcal{M}] E[\mathcal{M}^3] + 6 E^2[\mathcal{M}] \sigma^2[\mathcal{M}] + 3E^4[\mathcal{M}]}{\sigma^4[\mathcal{M}]} \\
\end{aligned}
\label{eq:pce03}
\end{equation}

In this work, we use the Weighted Least Squares (WLSQ) approach to compute the coefficients $c^i$. A set of $q$ samples of $\vksi$, \ $\vksi^{(i)}$ $(i=1,\dots,q)$, are defined and the value of the QoI at each sample is computed, forming  the vector $\mathbf{Q}=\left[\mathcal{M} \left(\vksi^{(1)}\right),\dots,\mathcal{M} \left(\vksi^{(q)}\right) \right]^\top \in \mathbb{R}^q$. The set of the unknown spectral coefficients $c^i$ forms the vector $\boldsymbol{c} = [c^0, \dots, c^P ]^{\top} \in \mathbb{R}^{P+1}$, and using this notation, the truncated expansion \eqref{eq:pce02} can be written as 
\begin{equation}
\mathbf{Q}  = \boldsymbol{\psi} \boldsymbol{c}  + \mathbf{\epsilon}, 
\label{eq:pce04}
\end{equation}
where $ \boldsymbol{\psi}$ is the measurement matrix with elements $\boldsymbol{\psi}_{ij} = \herm_j \left(\vksi^{(i)}\right) \in \mathbb{R}^{q \times (P+1)}$, i.e.\ the $i$-th row contains the values of the orthogonal polynomial bases $\herm_j$ evaluated at the $i$-th sample point, $\vksi^{(i)}$. In expanded form, $ \boldsymbol{\psi}$ is written as
\begin{equation}
\boldsymbol{\psi} = 
\begin{bmatrix}
\Psi_0 \left(\vksi^{(1)}\right)   &  \dots  & \Psi_P\left(\vksi^{(1)}\right)  \\
  \vdots     &  \ddots & \vdots    \\
\Psi_0\left(\vksi^{(q)}\right)  &  \dots   & \Psi_P\left(\vksi^{(q)}\right)   \\
\end{bmatrix}.
\label{eq:psi_matrix}
\end{equation}

The coefficient vector $\boldsymbol{c}$ is computed as a solution to the following weighted least-squares minimization problem, 
\begin{equation}
\min_{c} \left\lVert \boldsymbol{W}^{\frac{1}{2}} \Big ( \mathbf{Q} - \boldsymbol{\psi} \boldsymbol{c} \Big ) \right\rVert_2^2=\min_{c} \left ( \mathbf{Q} - \boldsymbol{\psi} \boldsymbol{c} \right )^\top  \boldsymbol{W} \left( \mathbf{Q} - \boldsymbol{\psi} \boldsymbol{c} \right ),
\label{eq:pce05}
\end{equation}
where $\| . \|_2$ is the Euclidean norm, $\boldsymbol{W}^{\frac{1}{2}}$ is a diagonal positive-definite matrix (to be defined later) and  $\boldsymbol{W} = (\boldsymbol{W}^{\frac{1}{2}})^{\top} \boldsymbol{W}^{\frac{1}{2}}$. The solution $\hat{\boldsymbol{c}}$ satisfies the normal equation 
\begin{equation}
\left( \boldsymbol{\psi}^{\top} \boldsymbol{W} \boldsymbol{\psi} \right) \hat{\boldsymbol{c}} =   \boldsymbol{\psi}^{\top} \boldsymbol{W} \mathbf{Q}.
\label{eq:pce06}
\end{equation}
The computational cost of forming and solving the above system is determined by the number of samples $q$. When random Monte Carlo sampling that follows the PDF is employed to obtain $\vksi^{(i)}$, $q \gg P+1$ is required for the system \eqref{eq:pce06} to be well-conditioned. For the close connection between the LSQ approach as $q \to \infty$ and the non-intrusive Galerkin projection with Gaussian quadrature, see chapter 3 of the book \cite{KnioBook}.  

We use the following weighting functions for Gaussian and uniform input distributions, 
\begin{equation}
\begin{aligned}
w^{\frac{1}{2}}_{H}(\vksi) &= exp(-\| \vksi^{(j)} \|_2^2/4)  \mbox{   (Gaussian distribution)},\\
w^{\frac{1}{2}}_{L}(\vksi) &= \prod_{k=1}^m \left(1-{\ksi_k^{(j)}}^2\right)^{1/4} \mbox{   (Uniform distribution)},
\end{aligned}
\label{eq:weights_coherence_optimal_sampling}
\end{equation}
where the subscripts 'H' and 'L' refer to Hermite and Legendre polynomials respectively. These analytic expressions are derived via an approach known as asymptotic sampling \cite{HADIGOL2018382,HAMPTON2015363,HAMPTON201573}.

System \eqref{eq:pce06} will be augmented with information on the sensitivity of $\mathcal{M} (\vksi^{(i)})$ with respect to $\vksi^{(i)}$; this is detailed in the next section. The sensitivities with respect to all the elements of $\vksi^{(i)}$ i.e.\ all random variables, can be evaluated very efficiently with  a single adjoint computation, as explained in section \ref{sec:adjoint basics}.

%%%%%%%%%%%%%%%%%%%%%%%%%%%%%%%%%%%%%%%%%%%%%%%%%%%%%%%%%%%%%%%%%%%%%%%
%%%%%%%%%%%%%%%%%%%%%%%%%%%%%%%%%%%%%%%%%%%%%%%%%%%%%%%%%%
\section{Sensitivity-enhanced gPC}
%%%%%%%%%%%%%%%%%%%%%%%%%%%%%%%%%%%%%%%%%%%%%%%%%%%%%%%%%%
\label{segPC basics}
Taking the gradient of the truncated expansion  \eqref{eq:pce02_truncated} with respect to the $k$-th random variable at the $j$-th sample point, we get
\begin{equation}
\frac {d \mathcal{M}}{d \ksi_k^{(j)}}  =\sum_{i=0}^{P} c^{i} \frac{\partial \herm_i}{ \partial \ksi_k^{(j)}} +\eta_k(\vksi^{(j)}) \quad (j=1,\dots,q).
\label{eq:sepce01}
\end{equation}
For each $k$, we have a block of $q$ equations that can be written in matrix form as, 
\begin{equation}
\frac{d {\mathbf{Q}}}{d \ksi_k}=
\frac{\partial \boldsymbol{\psi} }{\partial \ksi_k}
\boldsymbol{c}+\boldsymbol{\eta}_k(\vksi) \quad (k=1,\dots,m),
\label{eq:sepce02}
\end{equation}
where 
 $\frac{d {\mathbf{Q}}}{d \ksi_k} =\left [\frac{d \mathcal{M}}{d \ksi_k^{(1)} }, \dots, \frac{d \mathcal{M}}{d \ksi_k^{(q)} } \right]^\top$ and  
\begin{equation}
\frac{\partial \boldsymbol{\psi} }{\partial \ksi_k}=
\begin{bmatrix}
\frac{\partial \Psi_0 \left(\vksi^{(1)}\right) }{\partial \ksi_k^{(1)}}   &  \dots  & \frac{\partial \Psi_P \left(\vksi^{(1)}\right)}{\partial \ksi_k^{(1)}}   \\
  \vdots     &  \ddots & \vdots    \\
\frac{\partial \Psi_0 \left(\vksi^{(q)} \right)}{\partial \ksi_k^{(q)}}   &  \dots   & \frac{\partial \Psi_P \left(\vksi^{(q)} \right)}{\partial \ksi_k^{(q)}}
\end{bmatrix}.\\
\label{eq:nabla_psi_matrix}
\end{equation}

We now stack together one block of $q$ equations of the form \eqref{eq:pce04} (arising from the evaluation of $\mathcal{M}\left(\vksi^{(j)}\right)$) and $m$ blocks, of $q$ equations each, of the form \eqref{eq:sepce02} (arising from the evaluation of sensitivity of $\mathcal{M}\left(\vksi^{(j)}\right)$ at all sample points with respect to each uncertain variable). In total, there are $(1+m)\times q$ equations that can be put in matrix form as
\begin{equation}
\mathbf{G}  = \boldsymbol{\phi} \boldsymbol{c}  + \mathbf{\theta}, 
\label{eq:se_pce04}
\end{equation}
where
\begin{equation}
\begin{aligned}
& \mathbf{G}  =  \left [\mathbf{Q}, \frac{d {\mathbf{Q}}}{d \ksi_1}, \dots, \frac{d {\mathbf{Q}}}{d \ksi_m} \right ]^{\top}= \\
 & \left[\mathcal{M} \left(\vksi^{(1)}\right),\dots,\mathcal{M} \left(\vksi^{(q)}\right), \frac{d \mathcal{M}}{d \ksi_1^{(1)} }, \dots, \frac{d \mathcal{M}}{d \ksi_1^{(q)}}, \dots, 
\frac{d \mathcal{M}}{d \ksi_m^{(1)}},\dots,\frac{d \mathcal{M}}{d \ksi_m^{(q)}}  \right]^\top \in \mathbb{R}^{(1+m)q \times 1}
\end{aligned}
\label{eq:expanded_G}
\end{equation}
is a block column vector, and 
\begin{equation}
\boldsymbol{\phi}=
\begin{bmatrix}
\boldsymbol{\psi} \\  \frac{\partial \boldsymbol{\psi} }{\partial \ksi_1} \\ \vdots \\ \frac{\partial \boldsymbol{\psi} }{\partial \ksi_m} 
\end{bmatrix}=
 \begin{bmatrix}
\Psi_0 \left(\vksi^{(1)}\right)   &  \dots  & \Psi_P\left(\vksi^{(1)}\right)  \\
  \vdots     &  \ddots & \vdots    \\
\Psi_0\left(\vksi^{(q)}\right)  &  \dots   & \Psi_P\left(\vksi^{(q)}\right)   \\
\frac{\partial \Psi_0 \left(\vksi^{(1)}\right) }{\partial \ksi_1^{(1)}}   &  \dots  & \frac{\partial \Psi_P \left(\vksi^{(1)}\right)}{\partial \ksi_1^{(1)}}   \\
  \vdots     &  \ddots & \vdots    \\
\frac{\partial \Psi_0 \left(\vksi^{(q)} \right)}{\partial \ksi_1^{(q)}}   &  \dots   & \frac{\partial \Psi_P \left(\vksi^{(q)} \right)}{\partial \ksi_1^{(q)}} \\
  \vdots     &  \ddots & \vdots    \\
\frac{\partial \Psi_0 \left(\vksi^{(1)} \right)}{\partial \ksi_m^{(1)}}   &  \dots   & \frac{\partial \Psi_P \left(\vksi^{(1)} \right)}{\partial \ksi_m^{(1)}}   \\
 \vdots     &  \ddots & \vdots \\
\frac{\partial \Psi_0 \left(\vksi^{(q)} \right)}{\partial \ksi_m^{(q)}}   &  \dots   & \frac{\partial \Psi_P \left(\vksi^{(q)} \right)}{\partial \ksi_m^{(q)}}   
\end{bmatrix}
 \in \mathbb{R}^{(1+m)q \times (P+1)}
 \label{eq:expanded_phi}
\end{equation}
is a row block matrix.

We now  define the following weighted least squares minimisation problem  
\begin{equation}
\min_{c}  \| \boldsymbol{W'}^{\frac{1}{2}} \left(\mathbf{G}- \boldsymbol{\phi} \boldsymbol{c}\right) \|_2^2 =\min_{c}  \left(\mathbf{G}- \boldsymbol{\phi} \boldsymbol{c}\right)^\top  \boldsymbol{W}' \left(\mathbf{G}- \boldsymbol{\phi} \boldsymbol{c}\right),
\label{eq:sepce03}
\end{equation}
where $ \boldsymbol{W'}$ is a block diagonal weighting matrix, consisting of $1+m$ blocks of $ \boldsymbol{W}$ i.e.\   $ \boldsymbol{W'}=diag ( \underbrace{\boldsymbol{W}, \dots, \boldsymbol{W}}_{{1+m} \mbox { blocks}} ) $. The above problem can be written in expanded form as
\begin{equation}
\min_{c} \left ( \mathbf{Q} -  \boldsymbol{\psi} \boldsymbol{c}  \right )^{\top} \boldsymbol{W} \left (  \mathbf{Q} -  \boldsymbol{\psi} \boldsymbol{c}  \right ) + \sum_{k=1}^{m} \left (\frac{d {\mathbf{Q}}}{d \ksi_k}  - \frac{\partial \boldsymbol{\psi} }{\partial \ksi_k}  \boldsymbol{c} \right )^{\top} \boldsymbol{W} \left ( \frac{d {\mathbf{Q}}}{d \ksi_k}-\frac{\partial \boldsymbol{\psi} }{\partial \ksi_k} \boldsymbol{c} \right ),
\label{eq:sepce033}
\end{equation}
which is a generalisation of the minimisation problem \eqref{eq:pce05}. 
% The derivative of which with respect to $\boldsymbol{c}$ is given by

% \begin{equation}
% \frac{\partial \mathcal{F}}{\partial {\boldsymbol{c}}}= 2 \beta_1 \Big ( \boldsymbol{\psi} ^{\top}\boldsymbol{W} \boldsymbol{\psi} \boldsymbol{c} - \boldsymbol{\psi} ^{\top}\boldsymbol{W}\mathbf{Q}\Big ) + 2 \beta_2 \Big ( \nabla \boldsymbol{\psi} ^{\top}\boldsymbol{W} \nabla \boldsymbol{\psi} \boldsymbol{c}  - \nabla \boldsymbol{\psi} ^{\top}\boldsymbol{W} \nabla \mathbf{Q} \Big ).
% \label{eq:sepce04}
% \end{equation}
% \textcolor{red}{***sth is wrong, the rhs should contain $\boldsymbol{c}$****} \textcolor{green}{***you were right, i forgot the c, the expression is correct now. if you look at eq. (9) it is obvious where the c is missing, I doubled checked it from a least squares book to be sure****}

The vector  $\hat{\boldsymbol{c}}$ that minimises the above function is obtained from the solution of the normal set of equations  
\begin{equation}
\left( \boldsymbol{\phi}^{\top}   \boldsymbol{W'} \boldsymbol{\phi} \right) \hat{\boldsymbol{c}} =   \boldsymbol{\phi} ^{\top}  \boldsymbol{W'} \mathbf{G},
\label{eq:sepce10}
\end{equation}
which is analogous to system \eqref{eq:pce06}.

%%%%%%%%%%%%%%%%%%%%%%%%%%%%%%%%%%%%%%%%%%%%%%%%%%%%%%%%%%
\section{Efficient computation of sensitivities with an  adjoint formulation}
%%%%%%%%%%%%%%%%%%%%%%%%%%%%%%%%%%%%%%%%%%%%%%%%%%%%%%%%%%
\label{sec:adjoint basics}
We now consider the computation of the sensitivities $\frac {d \mathcal{M}}{d \ksi_k^{(j)}} \;  (k=1 \dots m)$ at the $j$-th sampling point. We assume systems that are described by a set of partial  differential equations that after spatial  discretisation take the form 
\begin{equation}
   \R(\bu, \vksi) = 0,
\label{eq:sePCE05}
\end{equation}
where $\bu \in \mathbb{R}^{N_{\bu}}$ is the state vector and $\R: \mathbb{R}^{N_{\bu}} \times  \mathbb{R}^{m} \to \mathbb{R}^{N_{\bu}}$ varies smoothly with $\bu$ and $\vksi$. For simplicity, we suppose that \eqref{eq:sePCE05} describes mathematically a steady problem, i.e.\ $\bu$ does not depend on time, but the analysis can be easily extended to unsteady problems as well. For example, in the context of Computational Fluid Dynamics (CFD), equation \eqref{eq:sePCE05} represents the discretised steady-state Navier-Stokes and continuity equations, the state vector $\bu$ contains the velocity and pressure at discrete points in the computational domain, and the stochastic vector $\vksi$ describes uncertainties in the geometry, boundary conditions etc. The QoI can therefore be written as 
\begin{equation}
\label{eq:sePCE04}
   \mathcal{M} = \mathcal{M}(\bu(\vksi), \vksi):  \mathbb{R}^{N_{\bu}} \times  \mathbb{R}^{m} \to \mathbb{R},
\end{equation}
\noindent where $\bu(\vksi)$ is the solution of \eqref{eq:sePCE05} for one realisation of $\vksi$. 

Using the adjoint formulation, the sensitivities $\frac {d \mathcal{M}}{d \ksi_k^{(j)}} \;  (k=1 \dots m)$ at each sampling point $j$ can be computed at a cost which is independent of the dimension $m$ of the vector $\vksi$. Below we outline the basic steps of the discrete adjoint approach, more details can be found in \cite{Jameson1988} and the review paper \cite{Luchini2014AdjointAnalysis}. It is also possible to derive a continuous adjoint formulation, but this requires to use the continuous form of the equations from which the discrete system  \eqref{eq:sePCE05} is derived. The discrete formulation suffices for the purposes of this section. We employ the continuous formulation in one of the test cases of section \ref{sec:High Dimensional Systems}. 

We start by defining the Lagrangian $\mathcal{L} : \mathbb{R}^{N_{\bu}} \times  \mathbb{R}^{m} \to \mathbb{R}$ 
\begin{equation}
\label{eq:sePCE06}
    \mathcal{L}(\bu(\vksi), \vksi) = \mathcal{M}(\bu(\vksi), \vksi) +  \bl^{\top} \R(\bu(\vksi), \vksi), 
\end{equation}
where $\bl \in \mathbb{R}^{N_{\bu}}$ are known as Lagrange multipliers (or adjoint variables), and $\bu(\vksi)$ is the solution of \eqref{eq:sePCE05} at the $j$-th sampling point, $\vksi^{(j)} = [ \ksi_1^{(j)}, \dots, \ksi_m^{(j)}]$. Differentiating $\mathcal{L}(\bu(\vksi), \vksi)$ with respect to $\vksi$ yields
\begin{equation}
 \frac{d \mathcal{L}}{d \vksi} = \left( \frac{\partial \mathcal{M}}{ \partial \bu} \right)^\top \frac{\partial \bu}{ \partial \vksi} +  \frac{\partial \mathcal{M}}{ \partial \vksi} + \bl^{\top} \left ( \frac{\partial \R}{\partial \bu} \frac{\partial \bu}{ \partial \vksi} +  \frac{\partial \R}{\partial \vksi} \right ),
 \label{eq:sePCE07b}
\end{equation}
where $\frac{d \mathcal{L}}{d \vksi}=\left [ \frac{d \mathcal{L}}{d \ksi_1^{(j)}} \dots \frac{d \mathcal{L}}{d \ksi_m^{(j)}} \right ]^\top$ and $\frac{\partial \R}{\partial \bu}$ is the Jacobian matrix. The above expression can be re-arranged as,
\begin{equation}
 \frac{d \mathcal{L}}{d \vksi} = \left [ \left( \frac{\partial \mathcal{M}}{ \partial \bu}\right)^\top + \bl^\top \frac{\partial \R}{\partial \bu} \right ] \frac{\partial \bu}{ \partial \vksi} + \bl^{\top} \frac{\partial \R}{\partial \vksi}+ \frac{\partial \mathcal{M}}{ \partial \vksi}.
\label{eq:sePCE07}
\end{equation}

The Lagrange multipliers are chosen so that  $\frac{d \mathcal{L}}{d \vksi}$ becomes independent of variations of $\bu$ with respect to $\vksi$. This is achieved by setting the term within square brackets equal to $\mathbf{0}$, resulting in the Field Adjoint Equation (FAE) 
\begin{equation}
\label{eq:sePCE08}
\left( \frac{\partial \mathcal{M}}{ \partial \bu}\right)^\top + \bl^\top \frac{\partial \R}{\partial \bu}=\mathbf{0} \Longrightarrow 
    \left( \frac{\partial \R}{\partial \bu} \right)^\top \bl = -\frac{\partial \mathcal{M} }{\partial \bu}.
\end{equation}
Solving \eqref{eq:sePCE08}, the Lagrangian multiplies $\bl$ are obtained and the sensitivities can be computed from 
\begin{equation}
\label{eq:sePCE09}
\frac {d \mathcal{M}}{d \vksi}= \frac{d \mathcal{L}}{d \vksi} = \bl^{\top} \frac{\partial \R}{\partial \vksi} + \frac{\partial \mathcal{M}}{\partial \vksi},
\end{equation}
where $\frac {d \mathcal{M}}{d \vksi}=\left [ \frac{d \mathcal{M}}{d \ksi_1^{(j)}} \dots \frac{d \mathcal{M}}{d \ksi_m^{(j)}} \right ]^\top$.

The state and adjoint variables, $\bu$ and $\bl$ respectively, have the same dimension (equal to $N_u$) and the solution algorithms for the two sets of equations, \eqref{eq:sePCE05} and \eqref{eq:sePCE08}, are similar.  In the following, we will assume that the computational cost is also the same. Note that in practise the cost of solving \eqref{eq:sePCE08} is lower compared to that of \eqref{eq:sePCE05}, because the former set is linear, so no outer iterations to deal with the non-linearity are necessary. The cost of \eqref{eq:sePCE09} is considered negligible.

To summarise, the solution of equation \eqref{eq:sePCE08} at each sampling point provides a set of $m$ extra conditions at the cost of one additional evaluation. It is exactly this property of the adjoint formulation that is exploited in the se-gPC method to enrich efficiently the least squares system \eqref{eq:pce06}. 

For the lowest chaos order approximation with $p = 1$, the number of coefficients, see equation \eqref{theory03}, is equal to $P + 1= \frac{(1+m)!}{1!m!} = 1+m$. One direct and one adjoint evaluation can yield $1+m$ equations at a single sampling point, thus making the computational cost of the se-gPC independent from the dimension of the stochastic space $m$. 

For the standard least squares method, if $q$ sampling points are needed for a well conditioned least-squares system \eqref{eq:pce06}, $q$ equations are obtained at a cost of $q$ evaluations. For the same number of sampling points, the se-gPC will yield $(1+m)q$ equations, at a cost of $2q$ evaluations. If the number of equations is kept constant and equal to $q$, then se-gPC allows one to use a smaller number of sampling points, equal to $q/(1+m)$, with cost  $2q/(1+m)$ evaluations. For cases with a large number of uncertain parameters, this is a very significant reduction in computational cost. 

If the computational cost is kept constant (for example due to a fixed computational budget) and equal to $q$ evaluations, the number of sampling points with se-gPC will be $q/2$ and the number of equations $(1+m)q/2$. 

In all the above scalings, we have assumed that for each sampling point, both the direct and the adjoint problems are  solved. This need not be the case, for example for some points only the direct problem is solved, while in other points both problems are solved (note that the adjoint problem requires the solution of the direct problem). 

The minimum number of equations needed to compute the spectral  coefficients is equal to $P+1$ (this corresponds to oversampling ratio equal to 1), thus $P+1$ and $(P+1)/(1+m)$ sampling points are required for the standard gPC and se-gPC respectively. Of course, one can use larger oversampling ratios. In the present paper we are using the smallest number of sampling points. 

\begin{figure}[!ht]
\centering
\includegraphics[scale=0.40, clip]{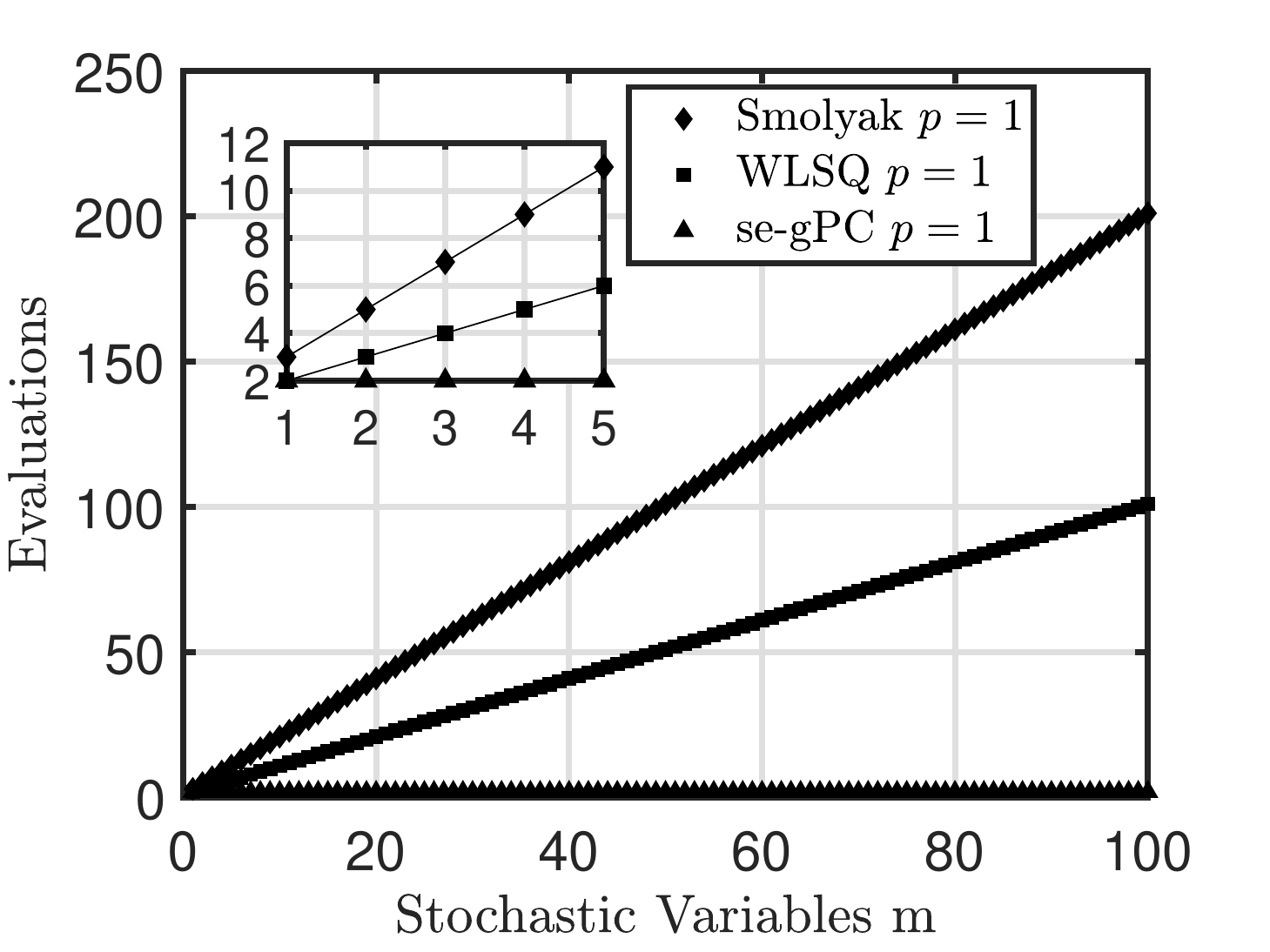}
\includegraphics[scale=0.40, clip]{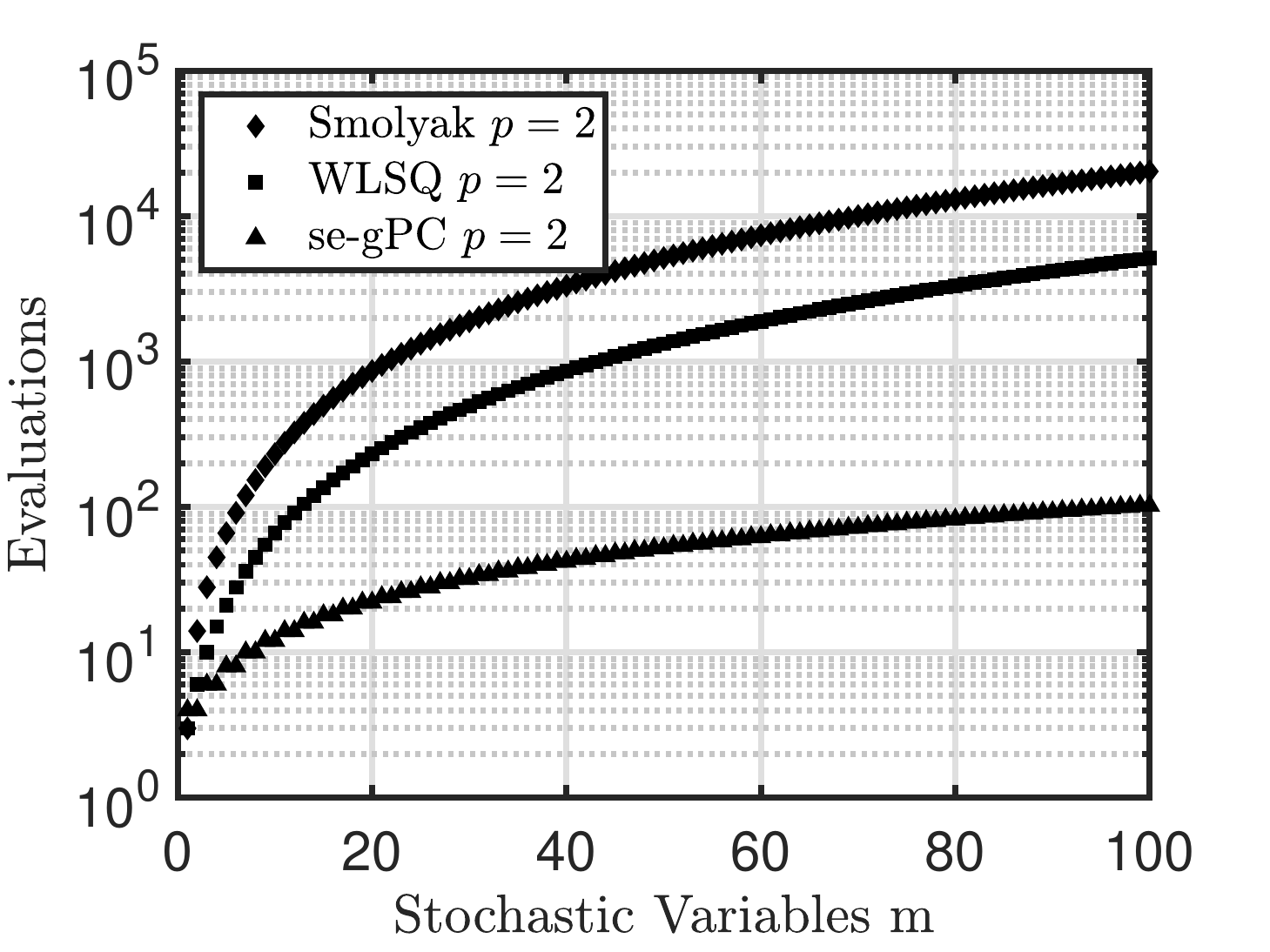}
\caption{Comparison of evaluations required by Smolyak grids, WLSQ and se-gPC at the QR points (see section \ref{sec:Fekete Sampling}), for $p = 1$ (left) and $p = 2$ (right).}
 \label{Tailored_Smol_Fek_se_comp}
\end{figure}

In figure \ref{Tailored_Smol_Fek_se_comp}, a comparison of the number of evaluations required by Smolyak as well as weighted Least Squares (WLSQ) and se-PCE (the latter two with the minimum number of sampling points) is made. Results are shown for two values of chaos order $p=1,2$. For $p = 1$ (left panel), the number of evaluations required by the se-gPC is constant and equal to 2, i.e.\ independent of the number of stochastic parameters, as mentioned earlier. For the other two methods, it scales  linearly with $m$; for WLSQ it is equal to $m+1$, and for Smolyak $2m+1$. For $p = 2$ (right panel), the number of evaluations scales  quadratically with $m$ for WLSQ and Smolyak (equal to $(m+1)(m+2)/2$ and $(m+1)(2m+1)$ respectively), but linearly for se-PCE (equal to $(m+2)$). This change in the scaling of the computational cost  with respect to $m$, from quadratic to linear, results in a significant  reduction (by a factor $m$) of the computational cost of se-PCE. This is precisely due to the property of the adjoints that provides $m$ conditions at the cost of just one evaluation. 

Since in the paper we are using the smallest number of sampling points, their selection needs to be done  carefully; this is examined in the next section. 

%%%%%%%%%%%%%%%%%%%%%%%%%%%%%%%%%%%%%%%%%%%%%%%%%%%%%%%%%%%%%%%%%%%%%%%
%%%%%%%%%%%%%%%%%%%%%%%%%%%%%%%%%%%%%%%%%%%%%%%%%%%%%%%%%%%%%%%%%%%%%%%
%%%%%%%%%%%%%%%%%%%%%%%%%%%%%%%%%%%%%%%%%%%%%%%%%%%%%%%%%%
\section{Sampling of the Stochastic Space}
%%%%%%%%%%%%%%%%%%%%%%%%%%%%%%%%%%%%%%%%%%%%%%%%%%%%%%%%%%
\label{sec:Fekete Sampling}
There are many different algorithms to sample the stochastic space for least squares PCE, for an overview see \cite{HADIGOL2018382}. In this paper we apply the pivoted QR decomposition algorithm, which has been used for sparse sensor placement for flow or image reconstruction \cite{BruntonPaper}, or to find optimal evaluation points and weights for integration \cite{SOMMARIVA20091324}. This is a greedy algorithm, that maximises the determinant of an appropriate matrix; in this sense it is a D-optimal design method. For the application of this algorithm for sampling the stochastic space, see \cite{DIAZ2018640,BruntonPaper}. Below we summarise the central idea and implementation. 

Consider again the standard least-squares problem, equation \eqref{eq:pce04}, repeated below for convenience,
\begin{equation}
\mathbf{Q}  = \boldsymbol{\psi} {\boldsymbol{c}}  + \boldsymbol{\epsilon},
\label{eq:tailored03}
\end{equation}
and assume that the elements of the residual vector $\boldsymbol{\epsilon}$ are independent, zero mean Gaussian random variables with the same standard deviation, i.e.\ for the $i$-th element $\epsilon_i \sim N(0,\eta^2)$, with $i=1, \dots, q$. Assuming for the time being that the weighting matrix is identity, $\boldsymbol{W}=\boldsymbol{I}_q$, the solution $\hat{\boldsymbol{c}}$ of the least squares system satisfies $\left( \boldsymbol{\psi}^{\top} \boldsymbol{\psi} \right) \hat{\boldsymbol{c}} =   \boldsymbol{\psi}^{\top} \mathbf{Q}$, see \eqref{eq:pce06}. The covariance matrix of the difference  $\hat{\boldsymbol{c}}-\boldsymbol{c}$ is given by 
\begin{equation}
Var(\hat{\boldsymbol{c}}-\boldsymbol{c}) = \eta^2 \left ( \boldsymbol{\psi}^{\top} \boldsymbol{\psi} \right )^{-1} 
\end{equation}
This can be easily proved, 
\begin{equation}
\begin{aligned}
\hat{\boldsymbol{c}} - \boldsymbol{c}  & =  (\boldsymbol{\psi}^{\top}  \boldsymbol{\psi})^{-1}\boldsymbol{\psi}^{\top} \mathbf{Q}- \boldsymbol{c}  = (\boldsymbol{\psi}^{\top}  \boldsymbol{\psi})^{-1}\boldsymbol{\psi}^{\top} (\boldsymbol{\psi} \boldsymbol{c}  + \boldsymbol{\epsilon} )- \boldsymbol{c} \\
& = (\boldsymbol{\psi}^{\top}  \boldsymbol{\psi})^{-1}\boldsymbol{\psi}^{\top} \boldsymbol{\psi} \boldsymbol{c} + (\boldsymbol{\psi}  \boldsymbol{\psi})^{-1}\boldsymbol{\psi}^{\top}  \boldsymbol{\epsilon} - \boldsymbol{c} \\
& =  \boldsymbol{c}+ (\boldsymbol{\psi}^{\top}  \boldsymbol{\psi})^{-1}\boldsymbol{\psi}^{\top}  \boldsymbol{\epsilon} - \boldsymbol{c}  = (\boldsymbol{\psi}^{\top}  \boldsymbol{\psi})^{-1}\boldsymbol{\psi}^{\top}  \boldsymbol{\epsilon}
\end{aligned}
\label{eq:tailored04_proof}
\end{equation}
Setting $\boldsymbol{A} = (\boldsymbol{\psi}^{\top}  \boldsymbol{\psi})^{-1}\boldsymbol{\psi}^{\top} $, we get \begin{equation}
\begin{aligned}
& Var(\hat{\boldsymbol{c}} - \boldsymbol{c}) = \mathcal{E} \left[(\hat{\boldsymbol{c}} - \boldsymbol{c})(\hat{\boldsymbol{c}} - \boldsymbol{c})^\top \right]=\mathcal{E} \left[\left(\boldsymbol{A} \boldsymbol{\epsilon} \right) \left(\boldsymbol{A}\boldsymbol{\epsilon}\right)^\top\right]=\mathcal{E} \left [\boldsymbol{A}  (\boldsymbol{\epsilon}\boldsymbol{\epsilon}^\top) \boldsymbol{A} ^{\top} \right] = \\
& \boldsymbol{A}  \mathcal{E}(\boldsymbol{\epsilon} \boldsymbol{\epsilon}^{\top})\boldsymbol{A} ^{\top} 
= \boldsymbol{A}  \eta^2 \boldsymbol{I}_q \boldsymbol{A} ^{\top} =  \eta^2 \boldsymbol{A}   \boldsymbol{A} ^{\top} = \eta^2 (\boldsymbol{\psi}^{\top}  \boldsymbol{\psi})^{-1},
\end{aligned}
\label{eq:tailored04}
\end{equation}
\noindent where $\mathcal{E}(.)$ is the expectation operation. If we include the weighting matrix and follow the same steps, the covariance matrix becomes
\begin{equation}
Var(\hat{\boldsymbol{c}} - \boldsymbol{c}) = 
\eta^2 \left (  \boldsymbol{\psi}  ^{\top}  \boldsymbol{W}  \boldsymbol{\psi}  \right) ^{-1}.
\end{equation}

The measurement matrix $\boldsymbol{\psi}$ contains information for all $q$ sampling points, see \eqref{eq:psi_matrix}. The objective of the D-optimal design is to select the sampling points that minimise the determinant of the matrix $\left (\boldsymbol{\psi}  ^{\top}   \boldsymbol{\psi}  \right) ^{-1}$ or $\left (  \boldsymbol{\psi}  ^{\top}  \boldsymbol{W}  \boldsymbol{\psi}  \right) ^{-1}$. Here we follow a greedy approach, i.e.\ from a large pool of $q$ candidate sample points, we select the required number of points (equal to $P+1$ for the standard least squares method) that maximise the determinant. The process is as follows.

To select only a subset of the $q$ points, we multiply equation \eqref{eq:tailored03} with the row selection matrix $\boldsymbol{P} \in \mathbb{R}^{(P+1)\times q} $. At each row of $\boldsymbol{P} $ all elements are 0, except the element at the column that corresponds to the selected sampling point, which takes the value of 1. Thus \eqref{eq:tailored03} becomes, 
\begin{equation}
\mathbf{Q}_{\boldsymbol{P}} = \boldsymbol{P} \boldsymbol{\psi} {\boldsymbol{c}}  + \boldsymbol{\epsilon}_{\boldsymbol{P}},
\label{eq:tailored04b}
\end{equation}
where $\mathbf{Q}_{\boldsymbol{P}}=\boldsymbol{P}\mathbf{Q} $ now contains only the values of the QoI $\mathcal{M}$ at the selected $P+1$ sampling points. The covariance matrix is 
\begin{equation}
Var(\hat{\boldsymbol{c}} - \boldsymbol{c}) = \eta^2 \left [ \left ( \boldsymbol{P} \boldsymbol{\psi}\right) ^{\top} \left( \boldsymbol{P} \boldsymbol{\psi} \right ) \right]^{-1},    
\label{eq:tailored04a}
\end{equation}
and if we include the weighting matrix, it becomes
\begin{equation}
Var(\hat{\boldsymbol{c}} - \boldsymbol{c}) = 
\eta^2 \left [ \left ( \boldsymbol{P} \boldsymbol{W}^{1/2}  \boldsymbol{\psi} \right) ^{\top} \left( \boldsymbol{P}  \boldsymbol{W}^{1/2} \boldsymbol{\psi} \right ) \right]^{-1}
\end{equation}
The aim now is to chose the row-selection matrix $\boldsymbol{P}$ that minimizes the determinant of $\left [ \left ( \boldsymbol{P} \boldsymbol{W}^{1/2}  \boldsymbol{\psi} \right) ^{\top} \left( \boldsymbol{P}  \boldsymbol{W}^{1/2} \boldsymbol{\psi} \right ) \right]^{-1}$ or equivalently maximise the determinant of  the inverse, i.e.\
\begin{equation}
\pvt  =\underset{\boldsymbol{P}}{\mbox{argmax}} \mbox{ } det \left[ \left ( \boldsymbol{P} \boldsymbol{W}^{1/2}  \boldsymbol{\psi} \right) ^{\top} \left( \boldsymbol{P}  \boldsymbol{W}^{1/2} \boldsymbol{\psi} \right ) \right]
\label{eq:tailored05}
\end{equation}

The matrix $\boldsymbol{P}  \boldsymbol{W}^{1/2} \boldsymbol{\psi}\in \mathbb{R}^{(P+1)\times (P+1)}$ is square, thus
$det\left[ \left ( \boldsymbol{P} \boldsymbol{W}^{1/2}  \boldsymbol{\psi} \right) ^{\top} \left( \boldsymbol{P}  \boldsymbol{W}^{1/2} \boldsymbol{\psi} \right ) \right]= \left[ det \left( \boldsymbol{P}  \boldsymbol{W}^{1/2} \boldsymbol{\psi} \right )\right] ^2$, and the maximisation problem \eqref{eq:tailored05} becomes
\begin{equation}
\pvt =\underset{\boldsymbol{P}}{\mbox{argmax}} \mbox{ }  \left |  det \left(\boldsymbol{P}  \boldsymbol{W}^{1/2} \boldsymbol{\psi}\right)  \right | 
\label{eq:tailored06}
\end{equation}
The solution to the above maximisation problem is given by the pivoted QR decomposition
\begin{equation}
 \left( \boldsymbol{W}^{\frac{1}{2}} \boldsymbol{\psi} \right)^\top \pvt = \boldsymbol{Q} \boldsymbol{R}.
\label{eq:tailored01}
\end{equation}
There are very robust subroutines to perform this decomposition; in this paper we have used the 'qr' command of MATLAB.  The permutation matrix $\pvt$  is chosen so that the diagonal elements $R_{ii}$ of the upper diagonal matrix $\boldsymbol{R} $ are ordered in descending order, i.e. 
\begin{equation}
|R_{11}| \geq |R_{22}| \geq \dots \geq |R_{(P+1)(P+1)}|
\label{eq:tailored01a}
\end{equation}
The absolute value of the determinant is the product of the diagonal entries \cite{BruntonPaper},
\begin{equation}
\left | det \mbox{ }  \left( \boldsymbol{W}^{\frac{1}{2}} \boldsymbol{\psi} \right) \right | = \prod_{i=1}^{P+1} |R_{ii}|
\label{eq:tailored01_volume}
\end{equation}

In the present paper, we apply pivoted QR decomposition to the matrix $\left(\boldsymbol{W}^{\frac{1}{2}}\boldsymbol{\psi}\right)^\top$ of the standard weighted least-squares approach, not to the augmented matrix $\left( \boldsymbol{W'}^{\frac{1}{2}}\boldsymbol{\phi} \right) ^\top$  of the se-gPC approach; see \eqref{eq:expanded_phi} for the expanded form of $\boldsymbol{\phi}$. Of course, we could have applied it to the augmented matrix, but that would have identified not only the locations of the sampling points, but also what derivative of the QoI to take at each point in order to maximise the absolute value of the determinant of $\boldsymbol{P}\boldsymbol{W'}^{\frac{1}{2}}\boldsymbol{\phi} $. In the adjoint formulation however, all derivatives (sensitivities) are computed simultaneously, so the additional information provided by QR decomposition would not have been used.  Instead, we identify the $P+1$ points required by the standard LSQ approach, and then select the most important $(P+1)/(m+1)$ points (i.e.\ the ones with the highest diagonal elements $|R_{ii}|$) and apply on those the se-gPC.  

%%%%%%%%%%%%%%%%%%%%%%%%%%%%%%%%%%%%%%%%%%%%%%%%%%%%%%%%%%
\subsection{An example with two stochastic variables}
%%%%%%%%%%%%%%%%%%%%%%%%%%%%%%%%%%%%%%%%%%%%%%%%%%%%%%%%%%
%Weighted Tailored Sampling (WTS) produces a near-optimally conditioned measurement matrix $ \boldsymbol{\psi}^{\top} \boldsymbol{W} \boldsymbol{\psi}$.

\begin{figure}[!htb]
\centering
\begin{subfigure}[b]{0.43\textwidth} 	\includegraphics[scale=0.49, clip]{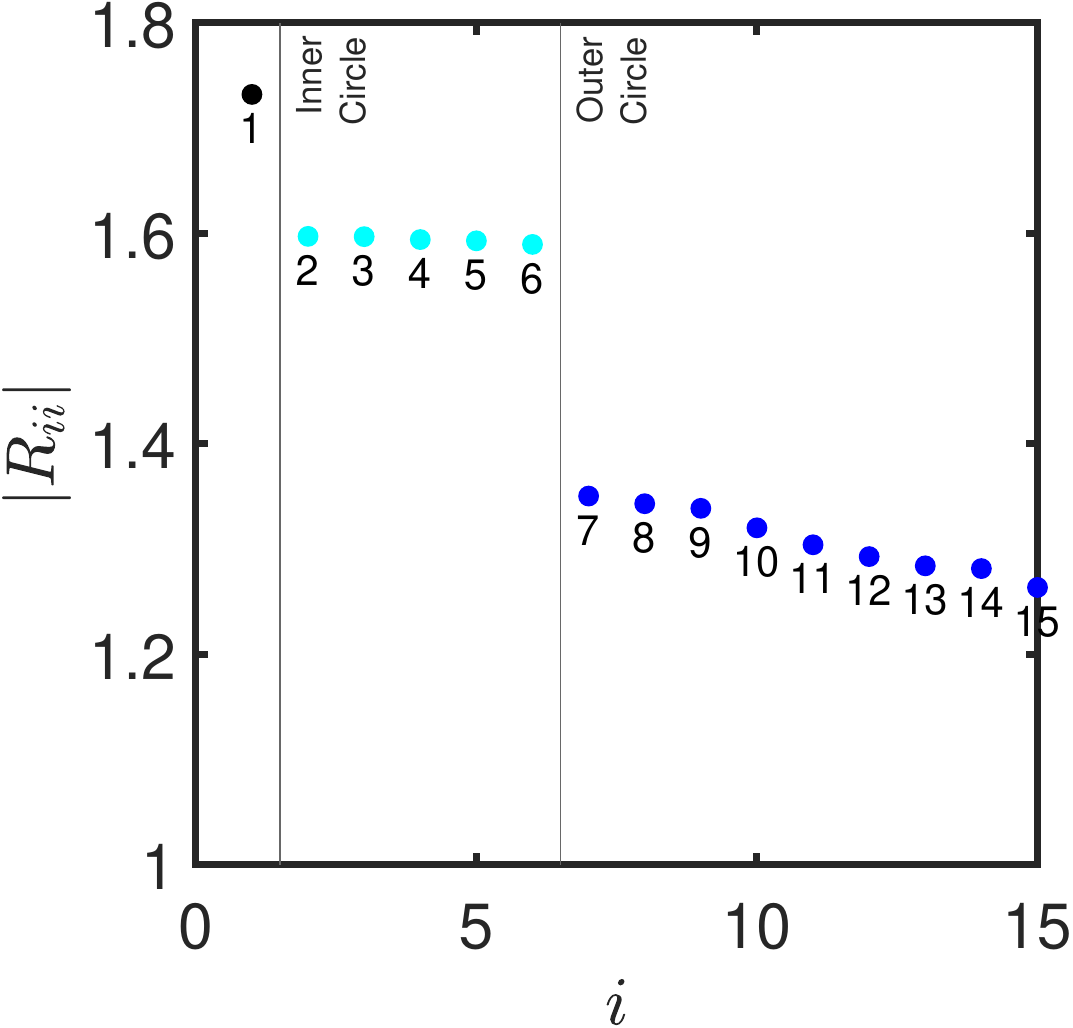}
\caption{}
\label{fig:grid_Rii_ranked}
\end{subfigure}
\begin{subfigure}[b]{0.43\textwidth}  \includegraphics[scale=0.49, clip]{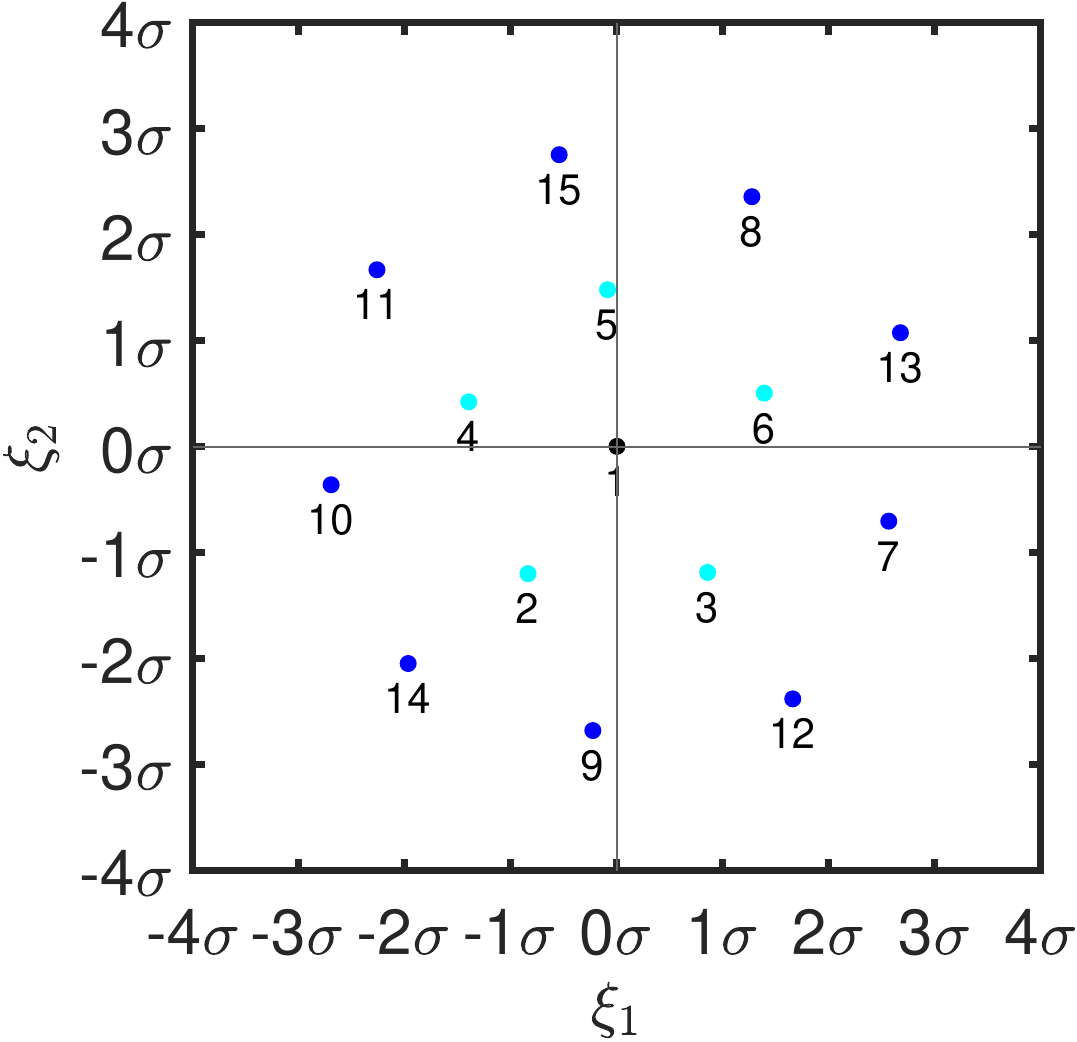}
\caption{}
\label{fig:grid_WTS_points}
\end{subfigure}
\caption{(a) Absolute values of the diagonal elements $|R_{ii}|$. (b) Selected QR sample points. Results for chaos order $p = 4$. The stochastic input is $\vksi \in \mathbb{R}^2$ with $\ksi_{1}, \ksi_{2} \sim N(0,1)$ sampled at $10000$ random points.}
\label{Grid Comparison}
\end{figure}

We have applied the QR decomposition algorithm to the case of $m = 2$ uncertain input variables, $\ksi_{1},\ksi_{2}$, both following standard normal distribution,  $N(0,1)$. We selected chaos order $p=4$, that results in $P+1=15$ unknowns. A pool of $10000$ random sample points were generated with Monte-Carlo from the input PDF. The absolute values of the diagonal elements of the matrix $\boldsymbol{R}$, $|R_{ii}|$, are shown in figure \ref{fig:grid_Rii_ranked}. The elements are tugged with their index number, $i=1 \dots 15 $. On the right panel, figure \ref{fig:grid_WTS_points}, the locations of the selected $15$ sampling points are shown, also tagged with the corresponding index $i$. The points are arranged in three concentric circles. The inner circle contains only one point; this point is located at the mean input, i.e.\ $\ksi_{1}^{(1)}=0$, $\ksi_{2}^{(1)}=0$, and corresponds to the highest ranked element $|R_{11}|$. The next circle contains 5 points, corresponding to elements $|R_{22}| \dots |R_{66}|$; the points are located at an average radius $\sim 1.47 \sigma$. The outer circle contains 8 points, corresponding to elements $|R_{77}| \dots |R_{15}|$, and has a radius of $\sim 2.77 \sigma$. 

\begin{figure}[!ht]
\centering
\includegraphics[scale=0.45, clip]{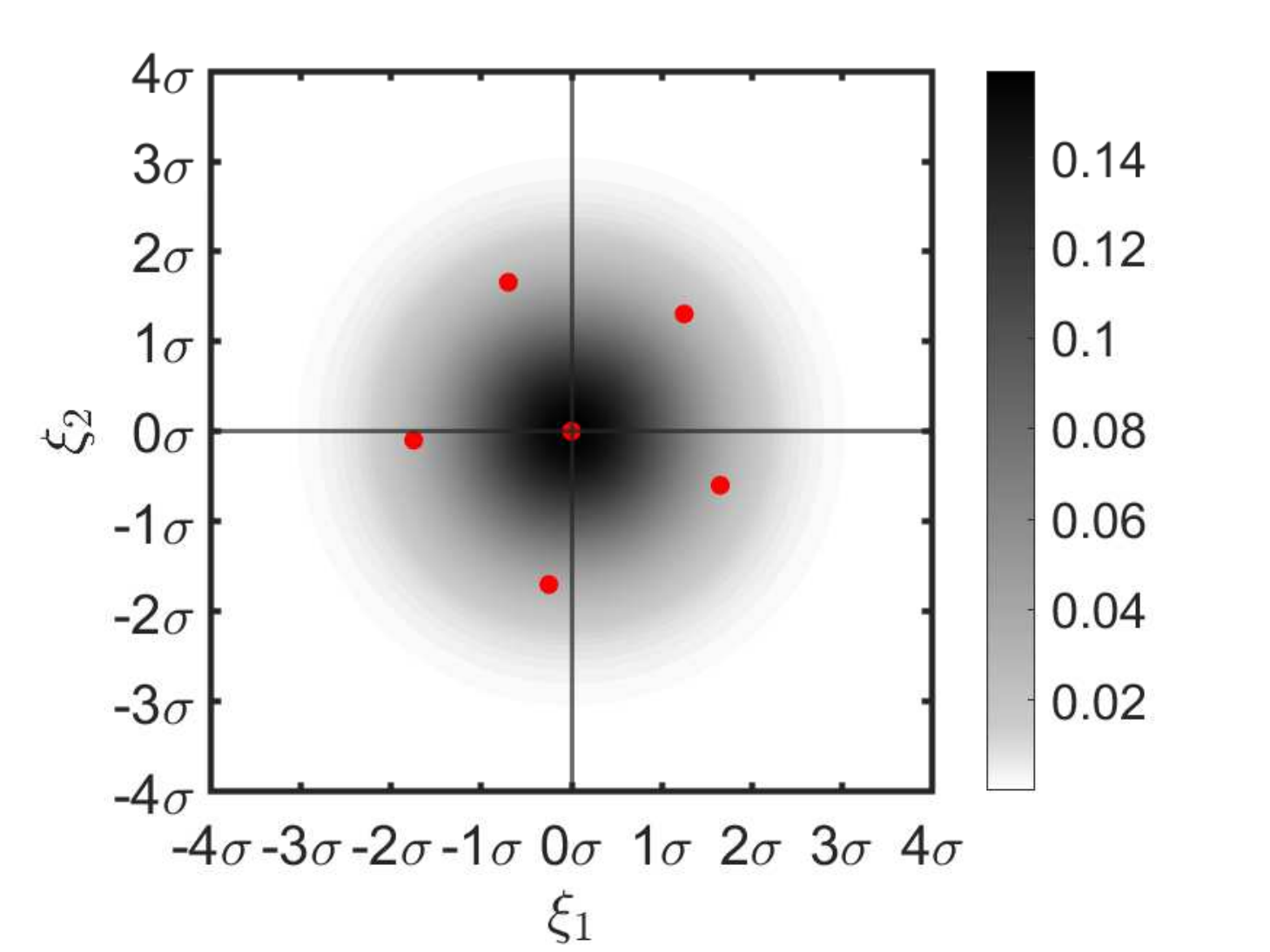}
\caption{Locations of QR sampling points for chaos order $p=2$, and two independent variables following standard  normal distribution, $\ksi_{1}, \ksi_{2} \sim N(0,1)$. The contour of the joint PDF is superimposed.  }
 \label{fig:WTS_TS_Comparison}
\end{figure}

We repeated the process for chaos order $p=2$ and the same input variables; now there are $P+1=6$ unknowns. The joint PDF and the selected QR sampling points are shown in figure \ref{fig:WTS_TS_Comparison}.  Again the first  point is placed at the peak of the PDF, i.e.\ $\ksi_{1}^{(1)}=0$, $\ksi_{2}^{(1)}=0$, and 5 points are symmetrically placed at a radius $\sim 1.75 \sigma$. 

We now proceed to compute the condition number (CN) of the matrix $\left(\boldsymbol{W}^{\frac{1}{2}}\boldsymbol{\psi}\right)^\top$. Note that the matrix appearing in  \eqref{eq:pce06} has CN which is the square of that of $\left(\boldsymbol{W}^{\frac{1}{2}}\boldsymbol{\psi}\right)^\top$. For $m=2$, and two chaos orders $p=1$ and $p=2$, all possible permutations of $P+1$ samples out of a pool of $q$ candidate samples (equal to $\frac{q!}{(q-P-1)!(P+1)!}$), were found and the $CN\left(\boldsymbol{W}^{\frac{1}{2}}\boldsymbol{\psi}\right)^\top$ was computed. We performed each numerical experiment 100 times; the averaged results are shown in table \ref{WTS_num_experiment}.  For example, for $p = 1$, $P+1 = 3$ and $q = 1000$, there are 166167000 permutations. The QR sampling points result in a system with condition number of $2.2$; only $0.24\%$ of all permutations have smaller CN. For $p=2$, the percentage is even smaller, $0.05\%$. 

\begin{table}[!ht]
\centering
\begin{tabular}{|l|l|l|l|l|}
\hline
$p$    & $P+1$ & $q$ & Number of permutations   & $CN \left(\boldsymbol{W}^{\frac{1}{2}}\boldsymbol{\psi}\right)^\top_{QR}$    \\ \hline \hline
$1$ & 3   & 1000 & 166167000    & 2.2 (0.24\%)      \\ \hline
$2$ & 6   & 50   & 15890700     & 16.5 (0.05\%) \\ \hline
\end{tabular}
\caption{Parameters for a numerical experiment to compute the CN of matrix $\left(\boldsymbol{W}^{\frac{1}{2}}\boldsymbol{\psi}\right)^\top$ evaluated at the QR sampling points.  The experiment was repeated 100 times and the average CN is reported in the last column. The percentage of $P+1$ sampling points from all permutations that have lower CN is also reported within brackets.}
\label{WTS_num_experiment}
\end{table}

\begin{figure}[!ht]
\centering
\begin{subfigure}[b]{0.49\textwidth} 	\includegraphics[scale=0.42, clip]{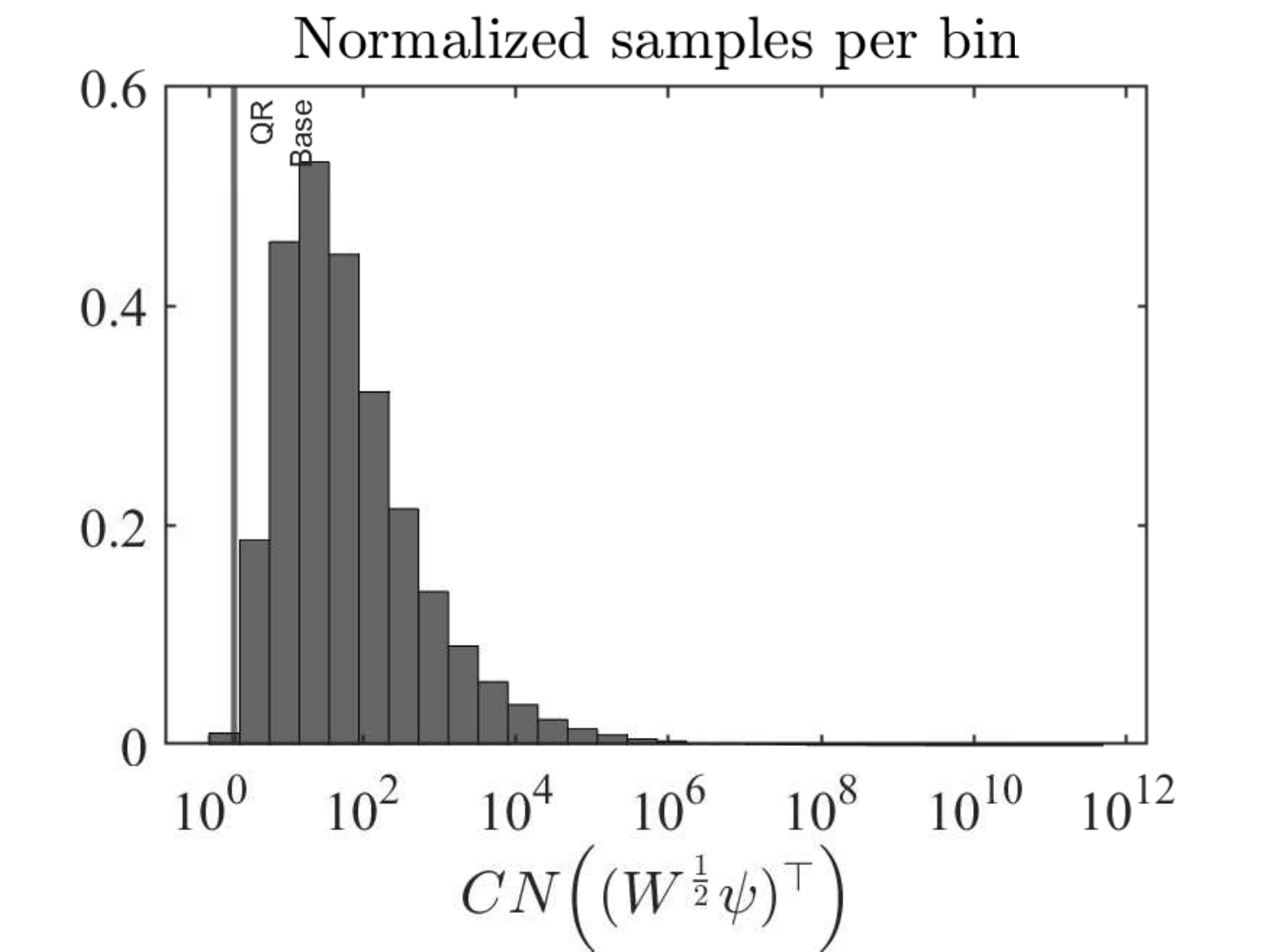}
\caption{}
\label{fig:num_exp_c1}
\end{subfigure}
\begin{subfigure}[b]{0.49\textwidth}  \includegraphics[scale=0.42, clip]{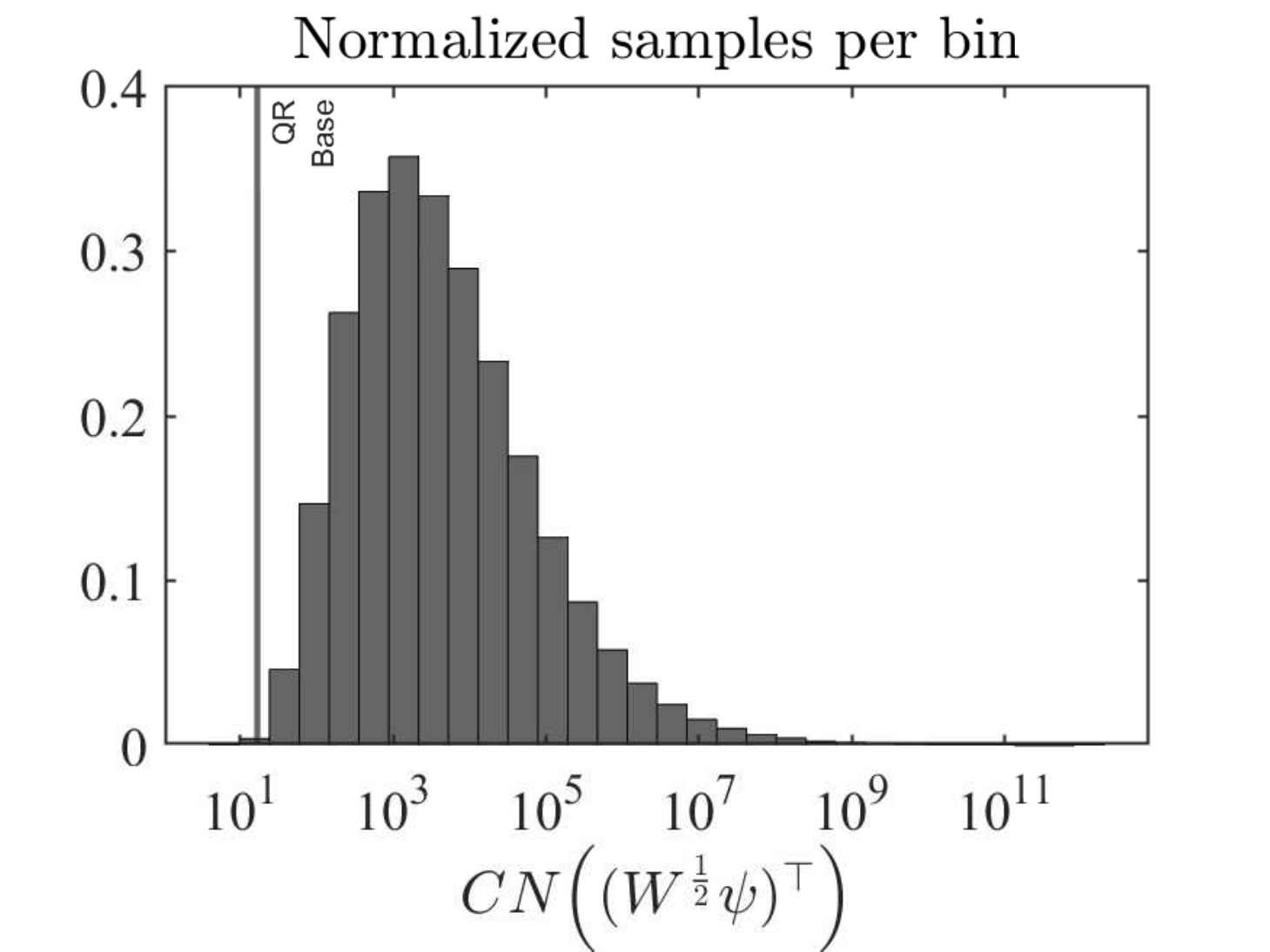}
\caption{}
\label{fig:num_exp_c2}
\end{subfigure}
\caption{Histogram of the CN (normalised samples per bin) for all permutations for $m = 2$ and (a) $p = 1$ and (b) $p = 2$. }
\label{Numerical Experiment}
\end{figure}

The histogram of the CN for all permutations is shown in figure \ref{Numerical Experiment}. Note that the horizontal axis is in logarithmic scale due to the very large spread of the CN values. The vertical lines mark the CN at the QR sampling points.  For $p=1$ (left panel), the median of the PDF is $\sim 43$, i.e. $20$ times larger. For $p=2$ (right panel) the PDF has similar shape, and the median CN is $2977$ i.e.\ 180 times higher than the one for the QR points.

In the following section, we apply standard WLSQ, se-PCE and Smolyak quadrature in several test cases and compare the three methods in terms of accuracy and computational cost.  

%%%%%%%%%%%%%%%%%%%%%%%%%%%%%%%%%%%%%%%%%%%%%%%%%%%%%%%%%%
\section{Application to test cases}
%%%%%%%%%%%%%%%%%%%%%%%%%%%%%%%%%%%%%%%%%%%%%%%%%%%%%%%%%%
\label{sec:High Dimensional Systems}
%%%%%%%%%%%%%%%%%%%%%%%%%%%%%%%%%%%%%%%%%%%%%%%%%%%%%%%%%%
\subsection{Stochastic Linear ODE}
%%%%%%%%%%%%%%%%%%%%%%%%%%%%%%%%%%%%%%%%%%%%%%%%%%%%%%%%%%
\label{ODE Linear Test Case}
A common benchmark for PCE methods is the linear ODE \cite{KarniadakisOriginal,GERRITSMA20108333},
\begin{equation}
\label{ODE_linear}
\frac{du}{dt} = -ku, \mbox {  } y(0) = 1,
\end{equation}
where $k$ is a single stochastic variable (so $m=1$), uniformly distributed in the interval $[0,1]$, i.e.\ the PDF is $f_k(k) = 1, \mbox{ } 0 \leq k \leq 1$, or $k \sim U(0,1)$.  For each $k$ realisation, the solution of \eqref{ODE_linear} is 
\begin{equation}
\label{ODE_linear Solution}
u(t;k) = e^{-kt},
\end{equation}
and the stochastic mean and variance can be computed analytically as,
\begin{equation}
\begin{aligned}
& \overline{u}(t) = \mathcal{E} \left [u(t;k) \right] = \int_0^1 e^{-kt} f_k(k) dk = \frac{1-e^{-t}}{t}\\
&Var \left(u(t;k)\right)= \mathcal{E} \left [\left(u(t;k)-\overline{u}(t)\right)^2 \right] = \int_0^1 \left ( e^{-kt} - \overline{u}(t)  \right )^2 f_k(k) dk = \frac{1-e^{-2t}}{2t}- \left ( \frac{1-e^{-t}}{t} \right )^2.
\end{aligned}
\label{eq:ODE_linear Stochastic Mean Variance}
\end{equation}
The quantity of interest is $u(t;k)$, and its derivative (sensitivity) with respect to $k$ is, 
\begin{equation}
\label{eq:ODE_linear Solution Derivative}
\frac{du(t; k)}{dk} = -te^{-kt}.
\end{equation}
For this test case, gPC provides accurate results for the mean and variance for small $t$, but as the time horizon $t$ grows, the results increasingly deviate from the analytical solution \eqref{eq:ODE_linear Stochastic Mean Variance}; for more details and the underlying fundamental physical reason, refer to \cite{GERRITSMA20108333}. 

We compute the spectral coefficients $c^i$ using Smolyak quadrature, as well as WLSQ and se-gPC at the QR sampling points. We consider a high chaos order $p=6$, and since $m=1$, there are in total $P+1=p+1=7$ coefficients. 

\begin{figure}[!ht]
\centering
\includegraphics[scale=0.40, clip]{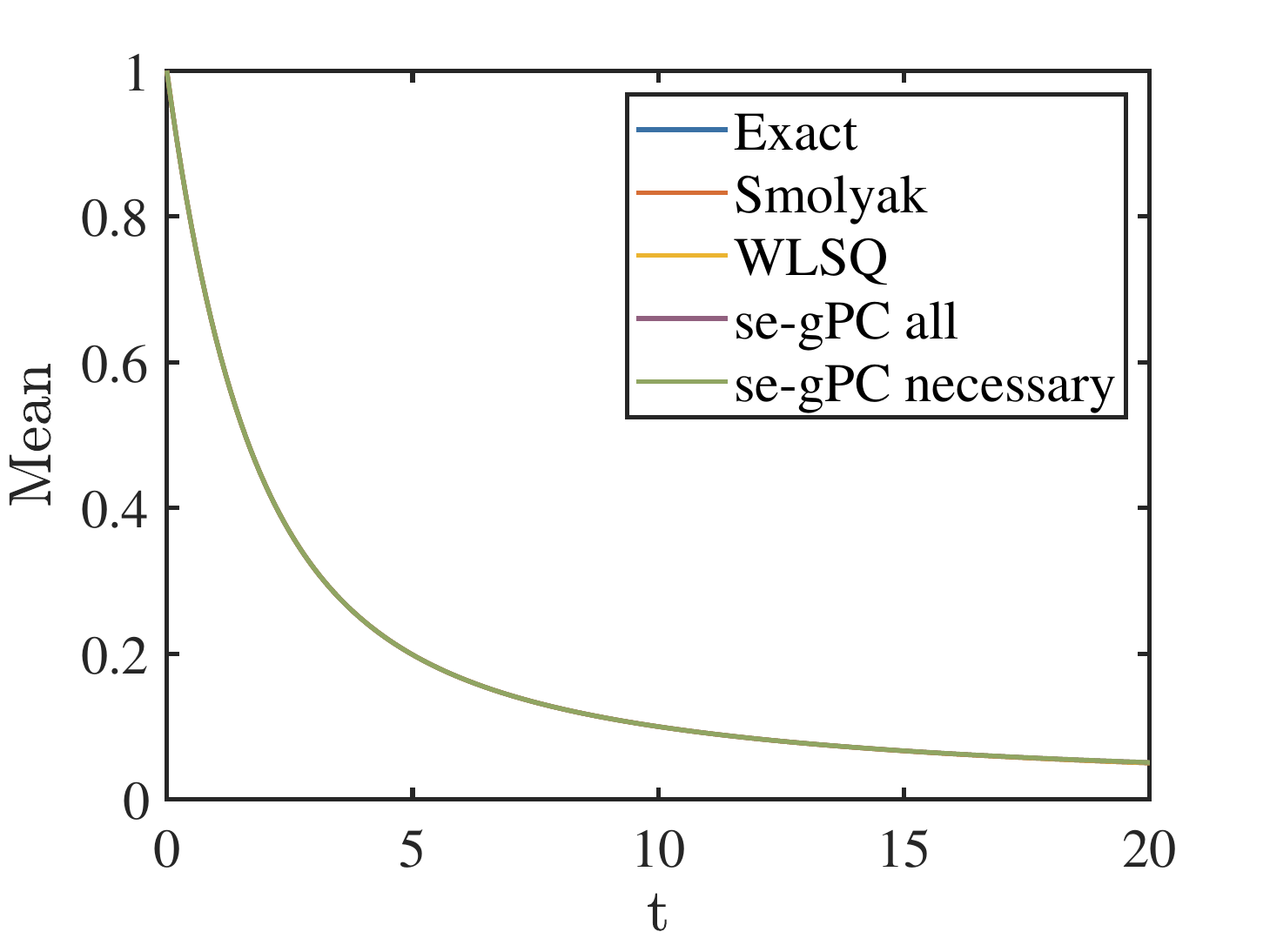}
\includegraphics[scale=0.40, clip]{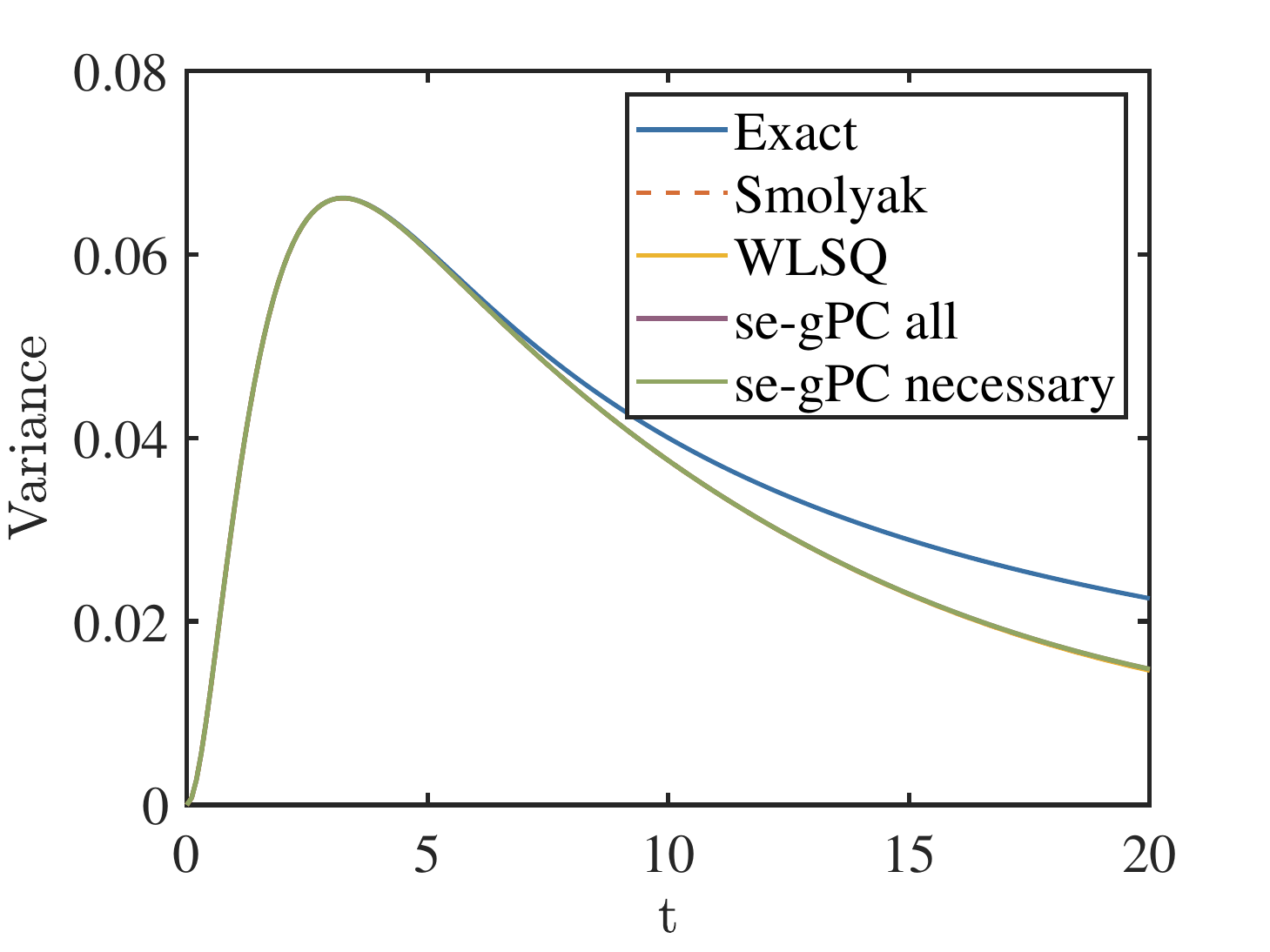}
\caption{Comparison between analytical solution (exact), Smolyak quadrature, WLSQ (at $7$ QR sampling points) and se-gPC  with sensitivity at $7$ QR sampling points (all), and at the top $4$ QR points (necessary) for $p = 6$.}
 \label{Exponential Comparison}
\end{figure}

The distribution of the mean and variance of $u(t;k)$ against time is shown in figure \ref{Exponential Comparison}. For se-gPC results from two cases are presented. In the first, all $7$ QR points were retained and additional $7$ sensitivity equations were added (this is the curve labelled "se-gPC all"). This method doubles the computational cost, i.e.\ requires 14 evaluations. For the second case, only $4$ QR points were utilised, at which both the QoI and its derivative were computed. The computational cost is $8$ evaluations (this is the curve labelled "se-gPC necessary"). This is the version of the se-gPC that will be used in the rest of the cases. Both versions produce results that match almost perfectly with standard WLSQ (at the 7 QR points) and Smolyak integration. This indicates that the se-gPC sampled at the minimum QR points can capture the statistics with an accuracy comparable with the other UQ approaches. Notice that for a small time window, say $t < 10$, the variance is predicted with high accuracy, but for larger $t$ the accuracy deteriorates, as mentioned earlier. 

The computational advantages of the method are not evident in this case because these is only one uncertain variable. Next, we consider cases with more uncertain variables.

%%%%%%%%%%%%%%%%%%%%%%%%%%%%%%%%%%%%%%%%%%%%%%%%%%%%%%%%%%
\subsection{ Ishigami Function}
%%%%%%%%%%%%%%%%%%%%%%%%%%%%%%%%%%%%%%%%%%%%%%%%%%%%%%%%%%
\label{Ishigami}
The Ishigami function \cite{Ishigami1990AnIQ},
\begin{equation}
\label{IshigamiFunction}
Y = sinX_1 + \alpha sin^2 X_2 + \beta X_3^4 sinX_1,
\end{equation}
is a commonly used benchmark case for UQ and Global Sensitivity Analysis \cite{SUDRET2008964,ABRAHAM2017461,ALLAIRE2012107,KERSAUDY2015103}. It has three stochastic variables $X_1$, $X_2$ and $X_3$ (so $m=3$) that follow a uniform distribution in the interval $[-\pi,\pi]$. This is a challenging case because the QoI, $Y$, is a strongly nonlinear and non-monotonic function, thus requires a high spectral order for accurate results; it is therefore an excellent case to assess the performance of se-gPC. 

\begin{figure}[!ht]
\centering
\includegraphics[scale=0.38, clip]{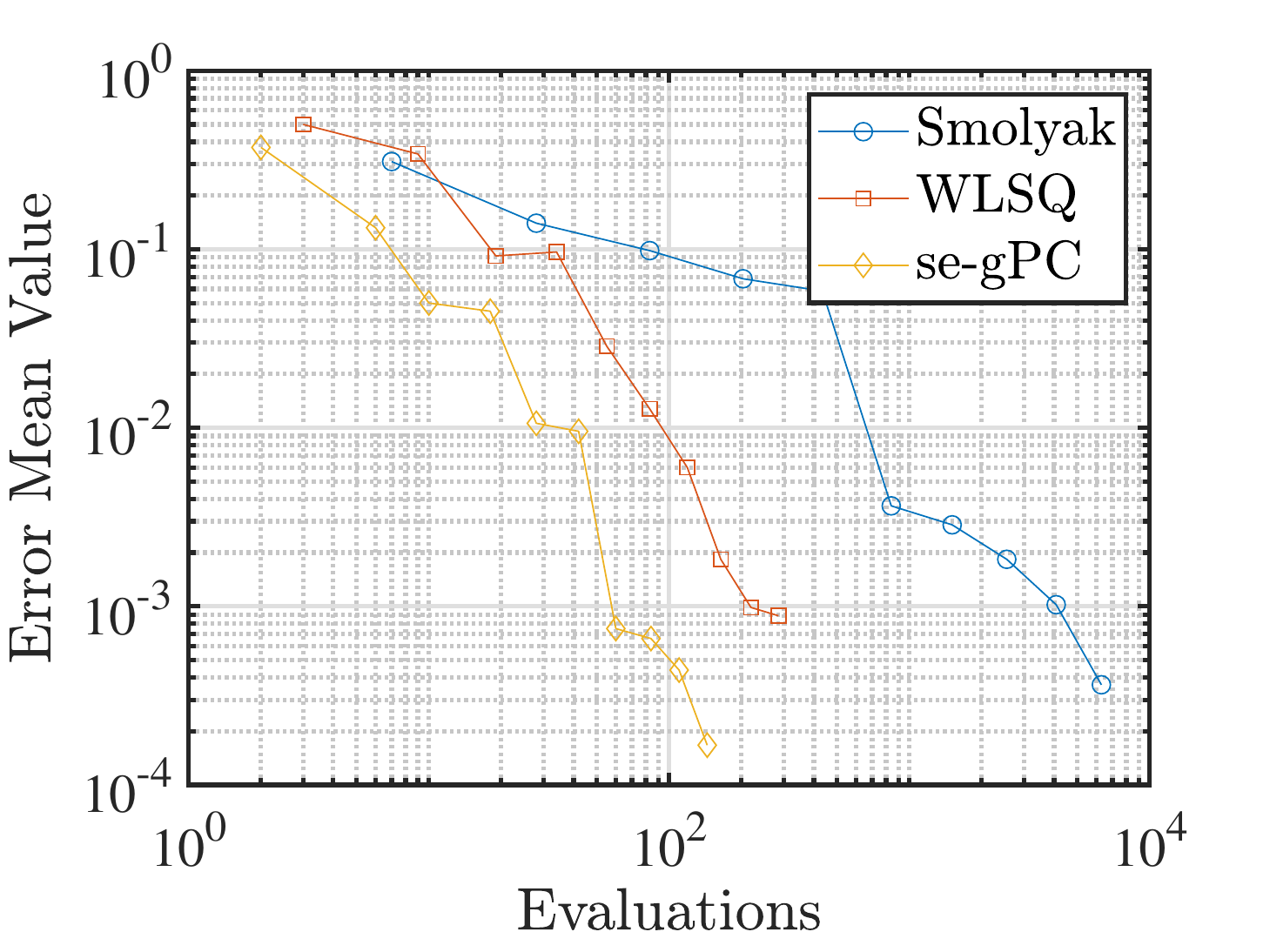}
\includegraphics[scale=0.38, clip]{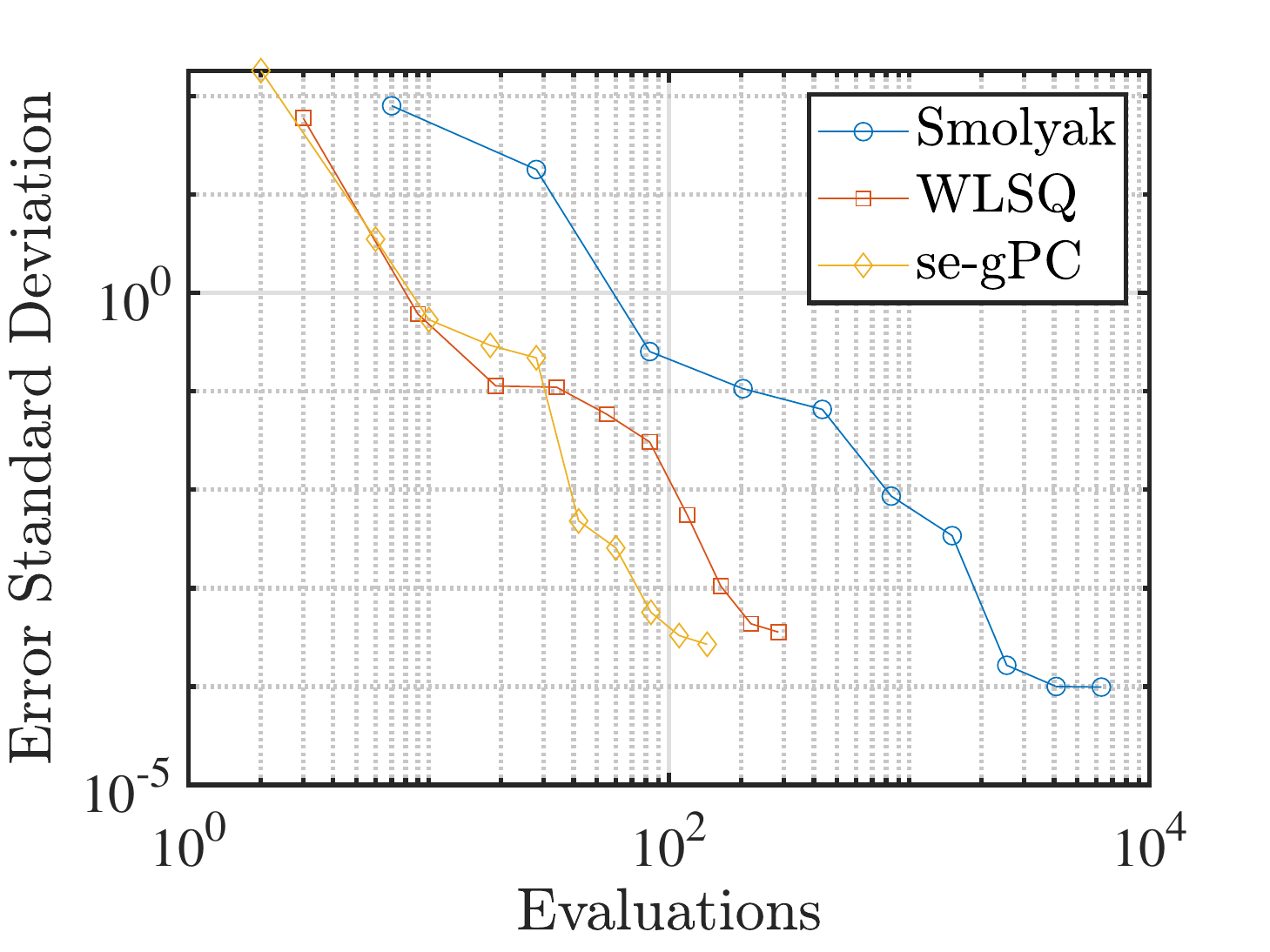}
\includegraphics[scale=0.38, clip]{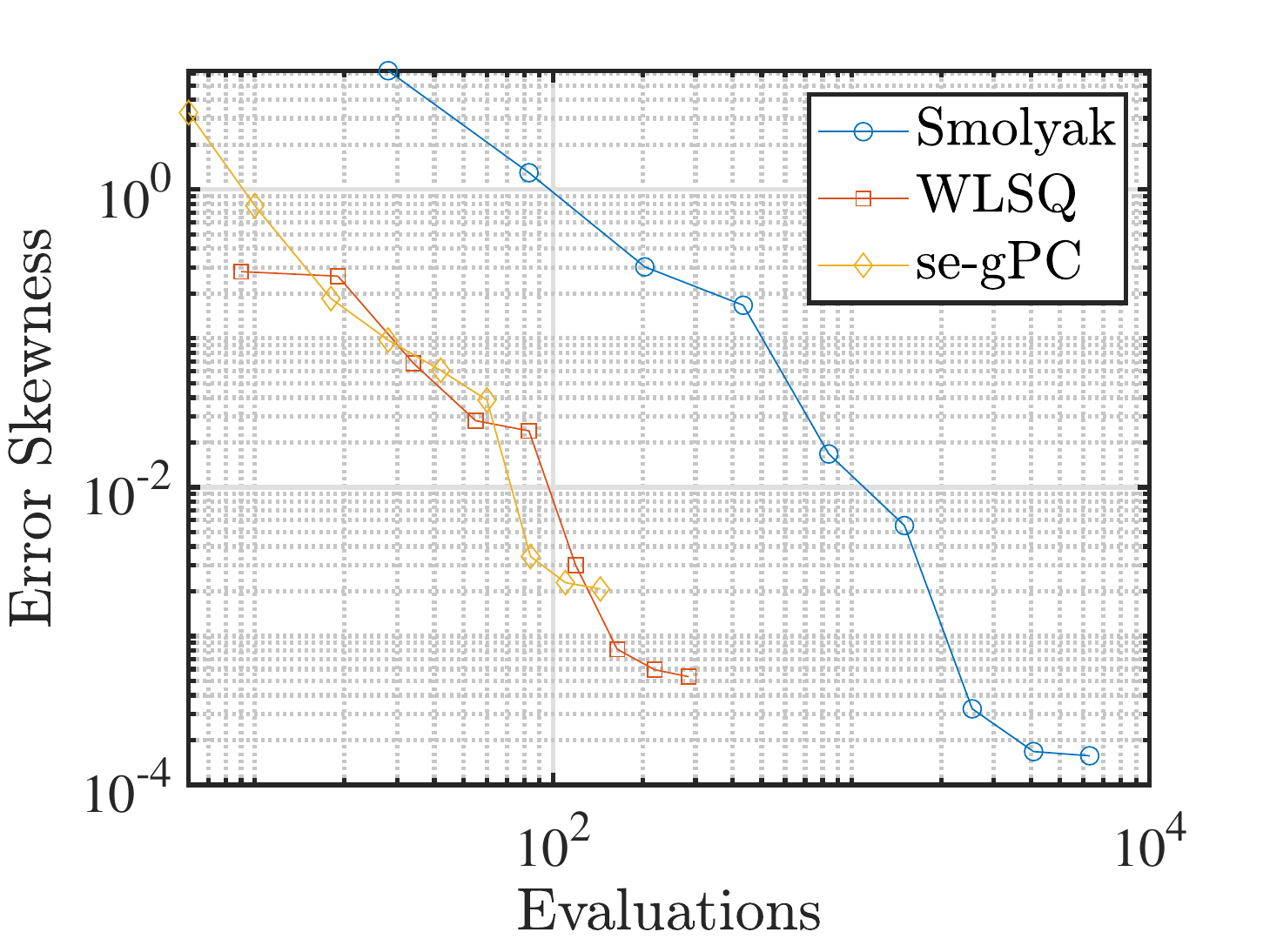}
\includegraphics[scale=0.38, clip]{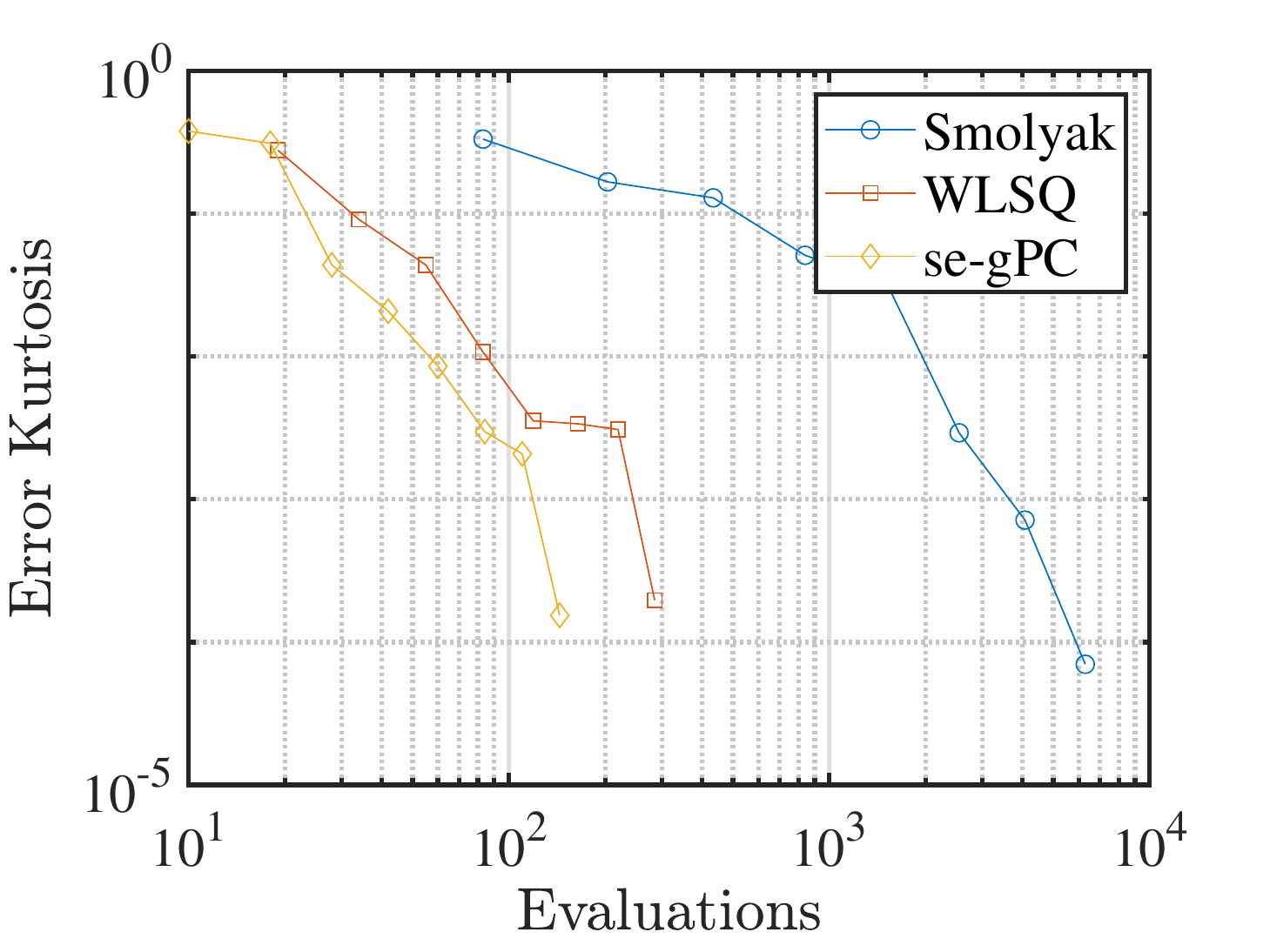}
\caption{ Convergence of the first four moments the Ishigami function (for $\alpha = 7$ and $\beta = 0.1$) against the number of evaluations, as computed with Smolyak Quadrature, WLSQ and se-GPC at the QR points. Each symbol represents a chaos order $p=1,\dots,10$.  Symbols are shown only for values of $p$ that allow the computation of the moments.}
 \label{Ishigami_Comparison}
\end{figure}

The statistical moments and sensitivity measures of $Y$ can be computed analytically. For $\alpha = 7$ and $\beta = 0.1$,  the reference values for the first four moments are $\widehat{\mu_Y} = 3.5000$, $\widehat{\sigma_Y} = 3.7208$, $\widehat{\delta_Y} = 0$ and $\widehat{\kappa_Y} = 3.5072$ \cite{BLATMAN2010183}. We computed the moments with Smolyak Quadrature as well as WLSQ and se-gPC at the minimum number of QR points and evaluated the error with respect to the reference values. For se-gPC, the sensitivities of QoI with respect to the 3 stochastic variables at the selected sampling points were computed analytically. The results are plotted against the number of evaluations in figure \ref{Ishigami_Comparison}. Each symbol in the plot represents a different chaos order $p=1, \dots,10$ (no symbols are shown for skewness and kurtosis for $p = 1$ and $p = 1,2$ respectively, because these moments cannot be computed for these values of $p$). As can be seen, all methods converge to very small errors. An interesting observation can be made from these plots. For small values of $p$, Smolyak quadrature can give slightly better results than se-gPC, see for example mean and standard deviation plots (top row). This is because Smolyak samples the QoI at more points, while for $p=1$ se-gPC samples the QoI and the sensitivities at just the first QR point. If high chaos order is needed for a challenging QoI like Ishigami,  information at a single point is not sufficient to produce accurate results. However, at larger $p$ values se-gPC samples the QoI and its sensitivities from more points, and this leads to rapidly improved convergence, as seen from the mean value and standard deviation plots. For the skewness and kurtosis, the superiority of se-gPC is not always evident compared to WLSQ, but both methods are more accurate compared to Smolyak quadrature for the same number of evaluations. Note that in this case, even though the derivatives of the Ishigami function are computed analytically at negligible  computational cost, the evaluations of the se-gPC are still counted as if an adjoint evaluation was performed (i.e.\ at each sampling point the cost is two evaluations). 

The Ishigami function is also commonly used for the evaluation of global sensitivity indices, like the Sobol indices \cite{BLATMAN2010183}, which quantify the effect of each input random variable onto the variance of the QoI. This information is useful for reducing the dimensionality of the input stochastic space by eliminating variables that have a small effect of the variance. The computed Sobol indices are shown in table \ref{IshigamiTable1} for two chaos orders; the se-gPC predictions match excellently with the analytic values for $p=6$. 

\begin{table}[!ht]
\begin{center}
\begin{tabular}{|l|l|l|l|}
\hline
Sobol index       &  Analytic value & se-gPC ($p = 6$)   & se-gPC ($p = 10$) \\ \hline \hline
$SU_1^{\top}$    &  $0.5574$ &     $0.5558$     &    $0.5567$   \\ \hline
$SU_2^{\top}$    &  $0.4424$ &     $0.4442$     &    $0.4433$      \\ \hline
$SU_3^{\top}$    &  $0.2436$ &     $0.2448$     &    $0.2449$      \\ \hline
\end{tabular}
\caption{Sobol indices comparison as computed by the gPC of the Ishigami function for $\alpha = 7$ and $\beta = 0.1$. The analytical values are obtained from \cite{BLATMAN2010183}.}
\label{IshigamiTable1}
\end{center}
\end{table}

%However, as discussed previously, the sensitivity enhanced gPC provides a more significant computational advantage to methods that have a lot of uncertain variables and can be approximated with relatively low chaos order. However it is useful to observe the convergence even for high chaos order. 

In the Ishigami function, where the number of uncertain variables is relatively small ($m=3$) and the required chaos order is high, the computational advantages of the se-gPC are not maximized, even though the method still outperforms the sparse quadrature. In the following, we consider cases with larger values of $m$.

%%%%%%%%%%%%%%%%%%%%%%%%%%%%%%%%%%%%%%%%%%%%%%%%%%%%%%%%%%
\subsection{ Two-dimensional viscous Burgers' equations}
%%%%%%%%%%%%%%%%%%%%%%%%%%%%%%%%%%%%%%%%%%%%%%%%%%%%%%%%%%
\label{Burgers2D}
We consider the two-dimensional viscous Burgers' equations, 
\begin{equation}
\label{Burgers2DPrimal}
\begin{aligned}
& u\frac{\pp u}{\pp x} + v\frac{\pp u}{\pp y} = \frac{1}{Re} \left ( \frac{\pp^2 u}{\pp x^2} + \frac{\pp^2 u}{\pp y^2} \right ),\\
& u\frac{\pp v}{\pp x} + v\frac{\pp v}{\pp y} = \frac{1}{Re} \left ( \frac{\pp^2 v}{\pp x^2} + \frac{\pp^2 v}{\pp y^2} \right ).
\end{aligned}
\end{equation}
This set is discretised in a square domain $\Omega = [0,1] \times [0,1]$ with Dirichlet boundary conditions at the top and bottom boundaries  $u(x,0) = u(x,1) = v(x,0) = v(x,1) = 0$ and Neumann boundary conditions, $\frac{\partial u}{\partial x}(1,y)=\frac{\partial v}{\partial x}(1,y)=0$ at the exit. At the inlet, the $u$ component of the velocity is a polynomial of degree $m+1$, i.e.\ $u(0,y) = \sum_{i=0}^{m+1} s_i y^i$, where the $m+2$ polynomial coefficients, $s_i$, are uncertain variables. Coefficients $s_1, s_2, \dots, s_m$ are free, while $s_0$ and $s_{m+1}$ are computed so that the boundary conditions at the top and bottom corners $u(0,0) = u(0,1) = 0$ are satisfied. The $v$ velocity at the inlet is given by $v(0,y)=-y^3 + y^2$, which satisfies the Dirichlet conditions $v(0,0) = v(0,1) = 0$. 

The QoI is the integral of the kinetic energy at the exit of the domain,
\begin{equation}
\label{Burgers2DObjective}
k_{e} =  \frac{1}{2} \int_0^1 \left( u^2 + v^2 \right) \Big |_{x=1} dy.
\end{equation}

The equations for the adjoint variables $u^{+}$ and $v^{+}$  can be easily derived and take the form,
\begin{equation}
\label{Burgers2DFAE}
\begin{aligned}
& u^{+} \frac{\pp v}{\pp y} + u\frac{\pp u^{+} }{\pp x} + v\frac{\pp u^{+} }{\pp y} +\frac{1}{Re} \left(\frac{\pp^2 u^{+}}{\pp x^2} + \frac{\pp^2 u^{+}}{\pp y^2} \right)  = v^{+} \frac{\pp v}{\pp x},\\
& v^{+} \frac{\pp u}{\pp x} + u\frac{\pp v^{+} }{\pp x} + v\frac{\pp v^{+} }{\pp y} +\frac{1}{Re} \left(\frac{\pp^2 v^{+}}{\pp x^2} + \frac{\pp^2 v^{+}}{\pp y^2} \right)  = u^{+} \frac{\pp u}{\pp y},
\end{aligned}
\end{equation} 
The boundary conditions for $u^{+}$ are  $u^{+}(0,y) = 0$, $u^{+}(x,0) = 0$ and $u^{+}(x,1) = 0$. Similarly for $v^{+}$, we get $v^{+}(0,y) = 0$, $v^{+}(x,0) = 0$ and $v^{+}(x,1) = 0$. %Note that the Dirichlet boundary conditions for $u$ and $v$ become Dirichlet boundary conditions for the adjoint variables $u^{+}$ and $v^{+}$.
The Neumann conditions at the exit, together with the form of the QoI, lead to Robin boundary conditions for the adjoint variables, 
\begin{equation}
\label{Burgers2DBCs_Neumann}
 u^{+}u + \frac{ 1}{ Re} \frac{\pp u^{+}}{\pp x} + u  = 0 \mbox{ and } v^{+}u + \frac{ 1}{ Re} \frac{\pp v^{+}}{\pp x} + v   = 0 , \mbox{ for } x = 1
\end{equation} 

The sensitivities of the QoI with respect to $s_i$  $(i=1,\dots,m)$ are given by
\begin{equation}
\label{Burgers2DSensitivities}
\frac{d k_{e}}{ds_i} = - \int_0^1 \left(u^{+}+v^{+}\right) u \Big |_{x=0}y^i dy - \frac{1}{Re} \int_0^1 \frac{\pp u^{+} }{\pp x}  \Bigg |_{x=0} y^i dy
\end{equation}

\begin{figure}[!ht]
\centering
\includegraphics[scale=0.4, clip]{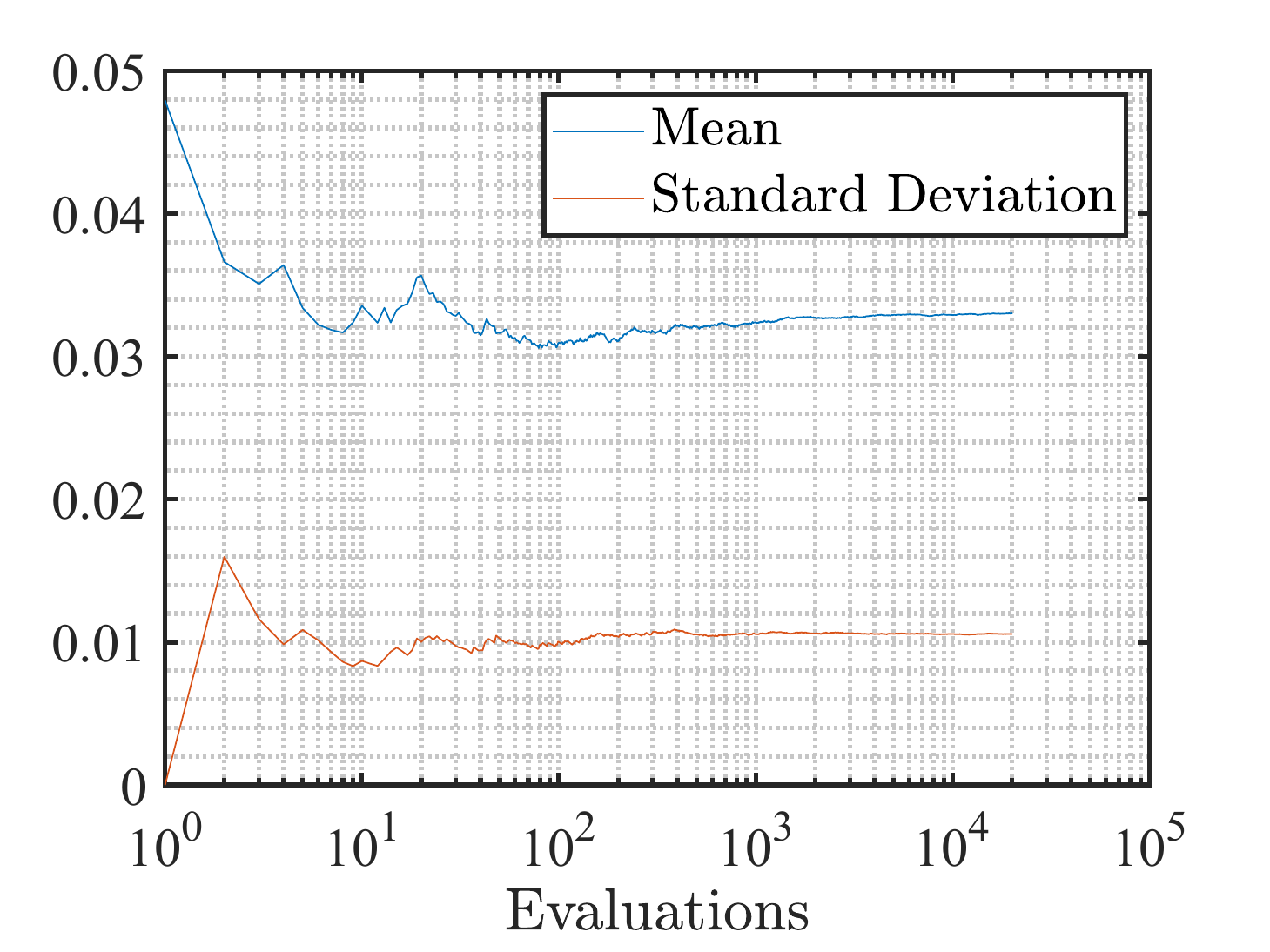}
\includegraphics[scale=0.4, clip]{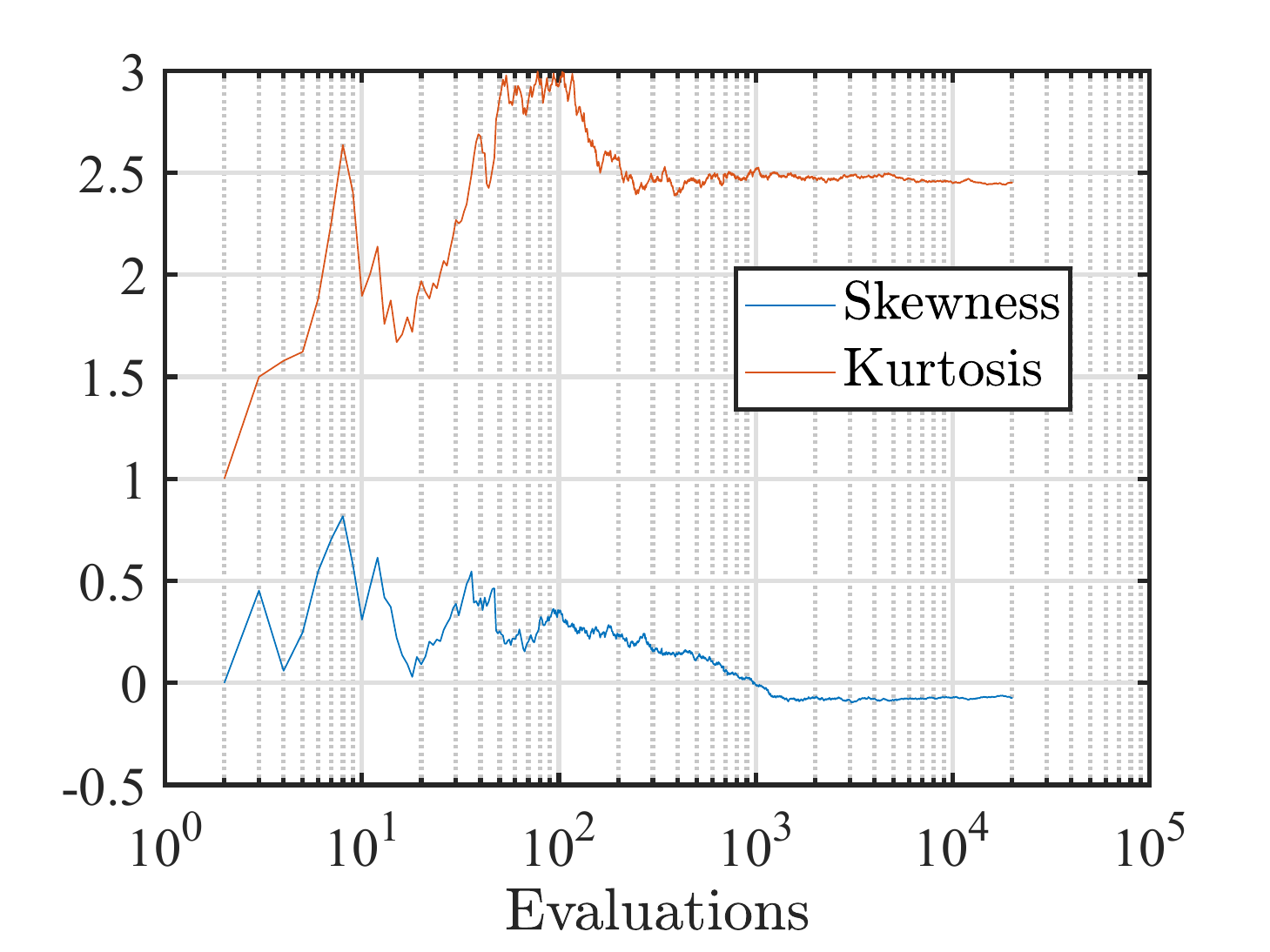}
\caption{Variation of the first four moments of $k_{\epsilon}$, see \eqref{Burgers2DObjective}, against the number of Monte-Carlo samples for $Re = 250$ and polynomial coefficients with mean $s_1 = -0.5$, $s_2 = -0.1$, $s_3 = 0.1$, $s_4 = 0.01$, $s_5 = -0.25$, $s_6 = 0.15$, $s_7 = 0.15$, $s_8 = -0.1$, $s_9 = 0.01$, $s_{10} = -0.25$, and standard deviation $\sigma_{s_i} = |s_i|/5$ ($i=1,\dots,10$). }
 \label{Burgers_Monte_Carlo_Convergence}
\end{figure}

We consider a case with $m = 10$ uncertain coefficients. The Reynolds number was set to $Re = 250$ and $\Omega$ was discretized using second-order finite differences  in a uniform $31 \times 31$ grid. To generate reference data for comparison, we performed $20,000$ Monte-Carlo simulations. The variation of the four moments of $k_{e}$ against the number of samples is shown in figure \ref{Burgers_Monte_Carlo_Convergence}.  The $3^{rd}$ and $4^{th}$ order moments, skewness and kurtosis respectively, have converged to one significant digit.  

\begin{figure}[!ht]
\centering
\includegraphics[scale=0.39, clip]{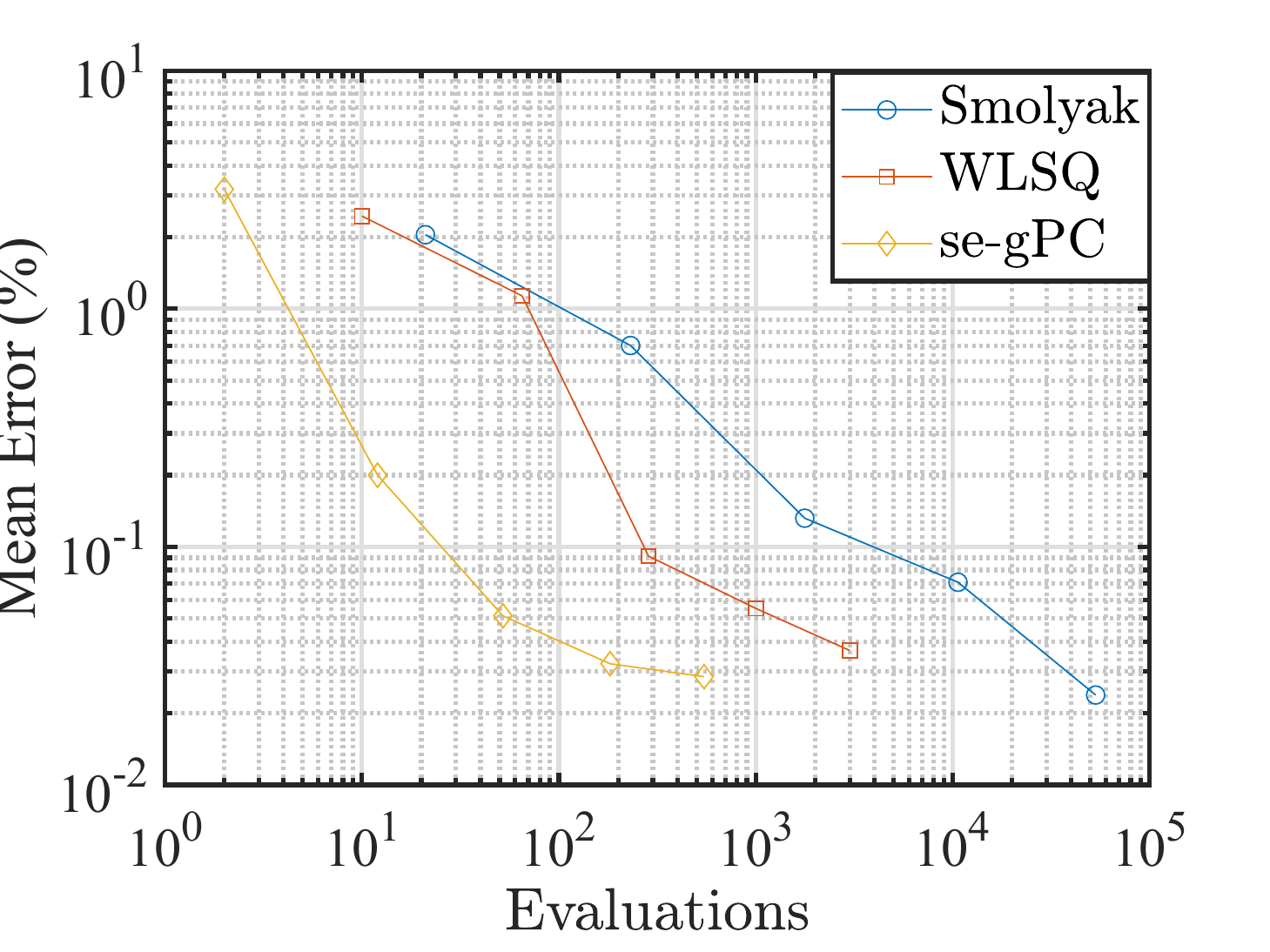}
\includegraphics[scale=0.39, clip]{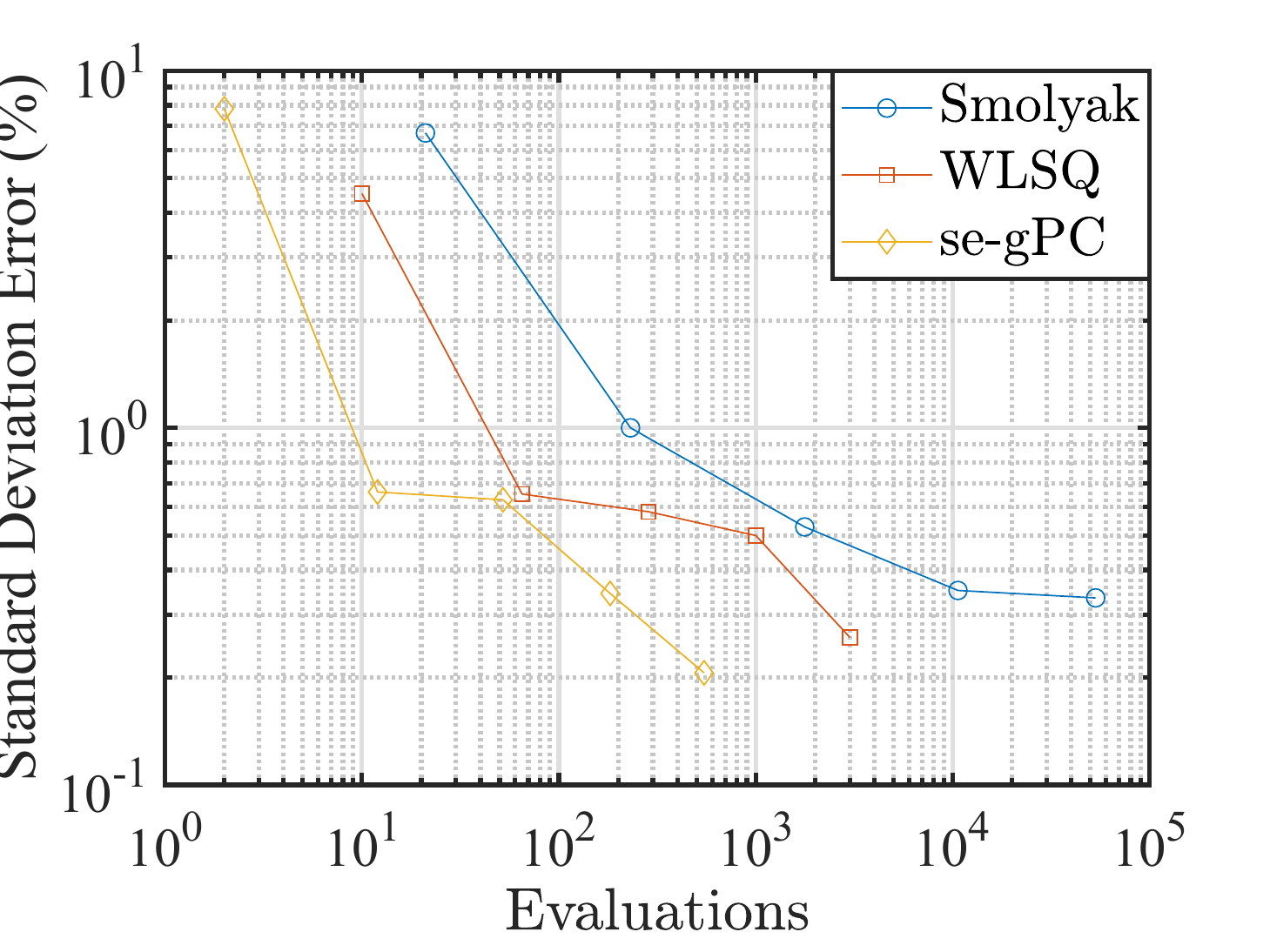}
\includegraphics[scale=0.39, clip]{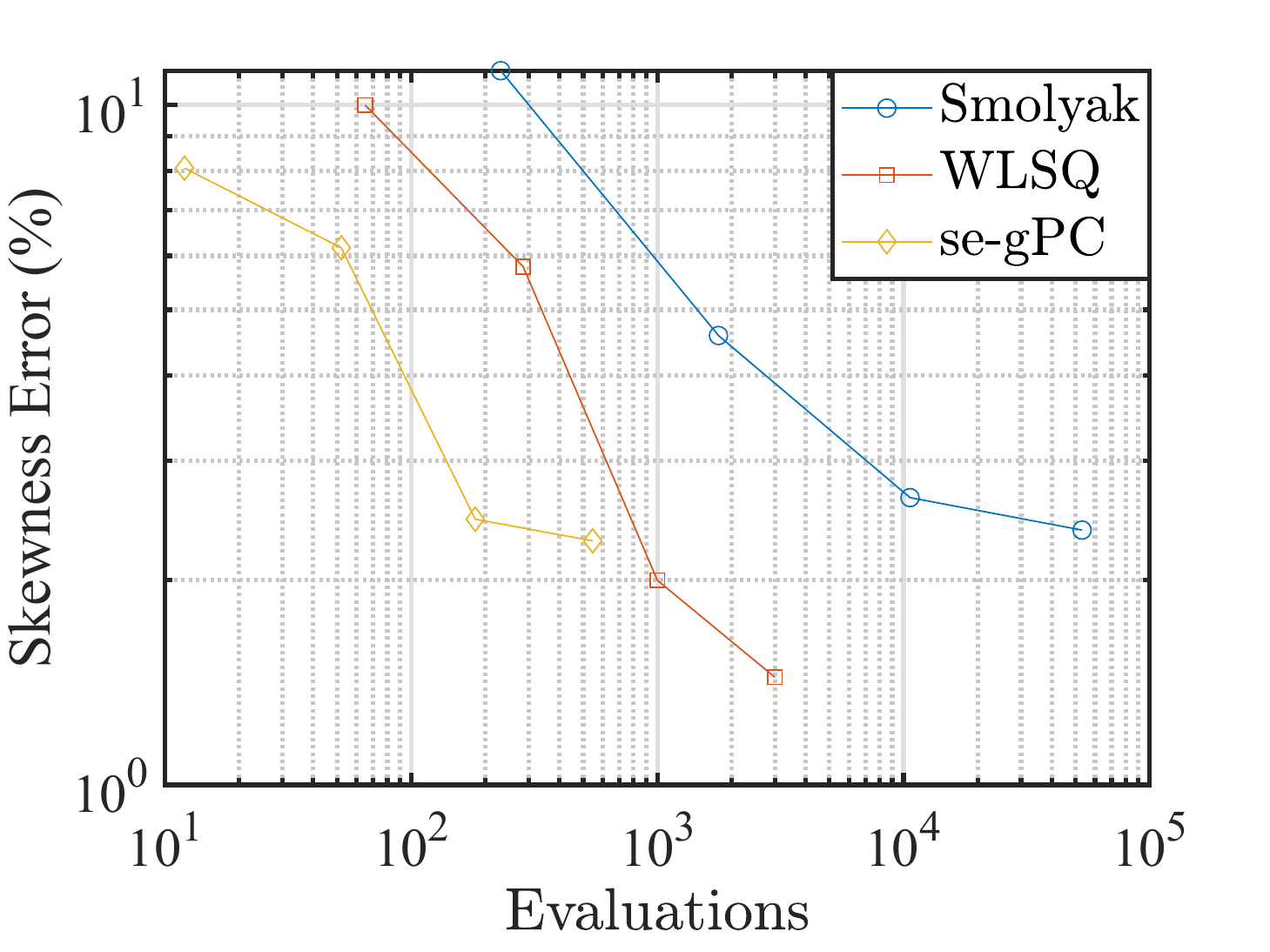}
\includegraphics[scale=0.39, clip]{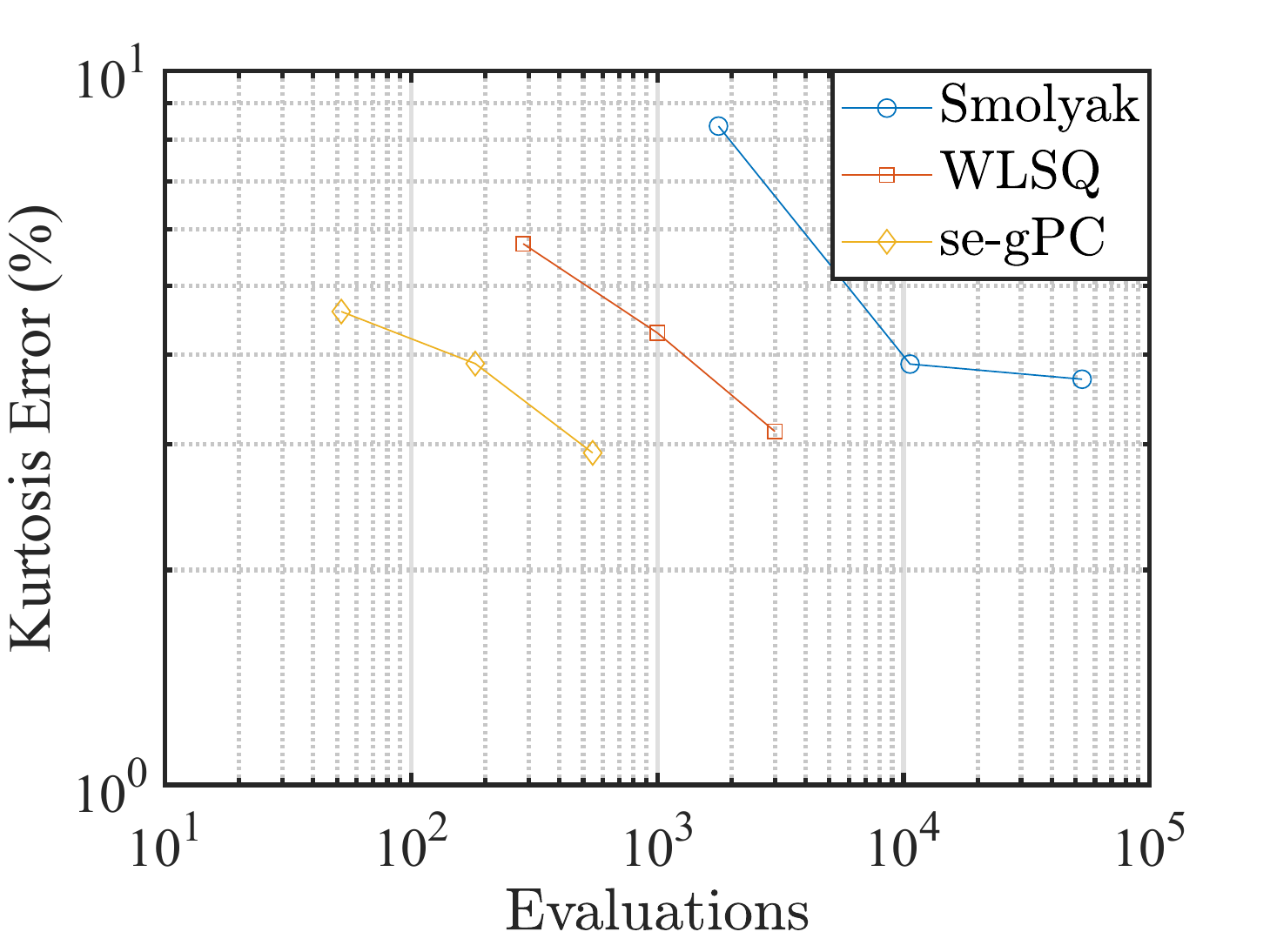}
\caption{Convergence of the first four moments of QoI $k_{\epsilon}$,  equation \eqref{Burgers2DObjective}, computed using Smolyak quadrature, WLSQ and se-gPC at the QR points, against number of evaluations. The points on each curve refer to a different chaos order. Symbols are shown only for values of $p$ that allow the computation of the moments.}
 \label{Burgers_Comparison}
\end{figure}

In figure \ref{Burgers_Comparison} the percentage  errors of the first four moments (with respect to the aforementioned Monte-Carlo statistics) as computed by Smolyak Quadrature, WLSQ and se-gPC with minimum QR points, are presented. Again each point represents a different chaos order. Symbols are shown only for values of $p$ that allow the computation of the moments. All methods converge, but the se-gPC achieves higher accuracy for the same computational cost compared to the other methods. For the same percentage error, the other two methods require one or two orders of magnitude more evaluations (note the logarithmic scale of the horizontal axis). Again it is possible that for low $p$, Smolyak may yield a lower error, but se-gPC converges faster as $p$ increases. The benefits of se-gPC are very significant for this case that has a relatively large number of stochastic variables. 

%%%%%%%%%%%%%%%%%%%%%%%%%%%%%%%%%%%%%%%%%%%%%%%%%%%%%%%%%%%%%%%%%%%%%%%

%%%%%%%%%%%%%%%%%%%%%%%%%%%%%%%%%%%%%%%%%%%%%%%%%%%%%%%%%%
\subsection{Transonic inviscid 2D flow around a NACA0012 airfoil}
\begin{figure}[!ht]
     \subfloat[\label{NACA12_Grid}]{\includegraphics[scale=0.3, clip]{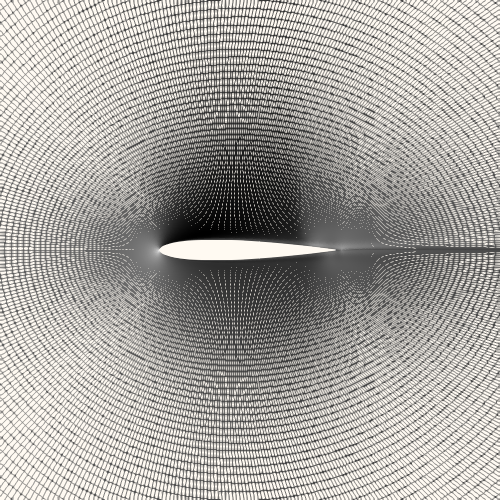}}
     \hfill
     \subfloat[\label{NACA12_flow}]{\includegraphics[scale=0.3, clip]{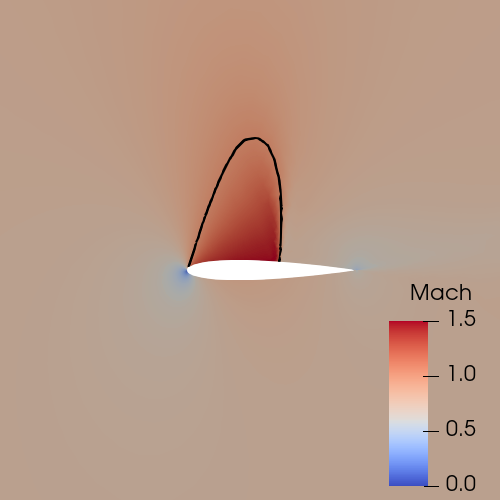}}
	\caption{Left: Unstructured mesh around the NACA0012 airfoil with $85000$ nodes (zoomed-in view close to the airfoil). Right: Contour plots of Mach number at $M_{\infty} = 0.75$ and $a_{\infty} = 4.0 ^{\circ}$. The black solid line indicates $M = 1$.}
	\label{MC_NACA2}
\end{figure}

%%%%%%%%%%%%%%%%%%%%%%%%%%%%%%%%%%%%%%%%%%%%%%%%%%%%%%%%%%
\label{Transonic Inviscid Naca}
This case considers the inviscid transonic $2D$ flow around a NACA0012 airfoil. The Euler solver of the SU2 package \cite{SU2} was used  to perform the numerical simulations for the direct and adjoint problems. The mesh consists of approximately $85,000$ nodes, it is refined close to the wall, and extends to $20$ chord lengths in the far field. A zoomed-in view is shown in figure \ref{NACA12_Grid}. At nominal free-stream conditions, the angle of attack is $a_{\infty} = 4.0 ^{\circ}$ and the Mach number is $M_{\infty} = 0.75$. For these conditions, the predicted lift coefficient is $C_L = 0.84$, in good agreement with experimental data \cite{SU2}. Contour plots of the Mach number are shown in figure \ref{NACA12_flow}. The flow is in the transonic regime, and a shock appears at the suction side of the airfoil. 

The QoI for this case is the lift coefficient, $C_L$. Two  scenarios are considered; uncertainty in the free-stream conditions, with or without uncertainty in the geometry.

%%%%%%%%%%%%%%%%%%%%%%%%%%%%%%%%%%%%%%%%%%%%%%%%%%%%%%%%%%
\subsubsection{ Uncertainty in the free-stream conditions}
%%%%%%%%%%%%%%%%%%%%%%%%%%%%%%%%%%%%%%%%%%%%%%%%%%%%%%%%%%
\label{Transonic Inviscid Naca Freestream Only}

\begin{figure}[!ht]
    \centering
    \includegraphics[scale=0.60, clip]{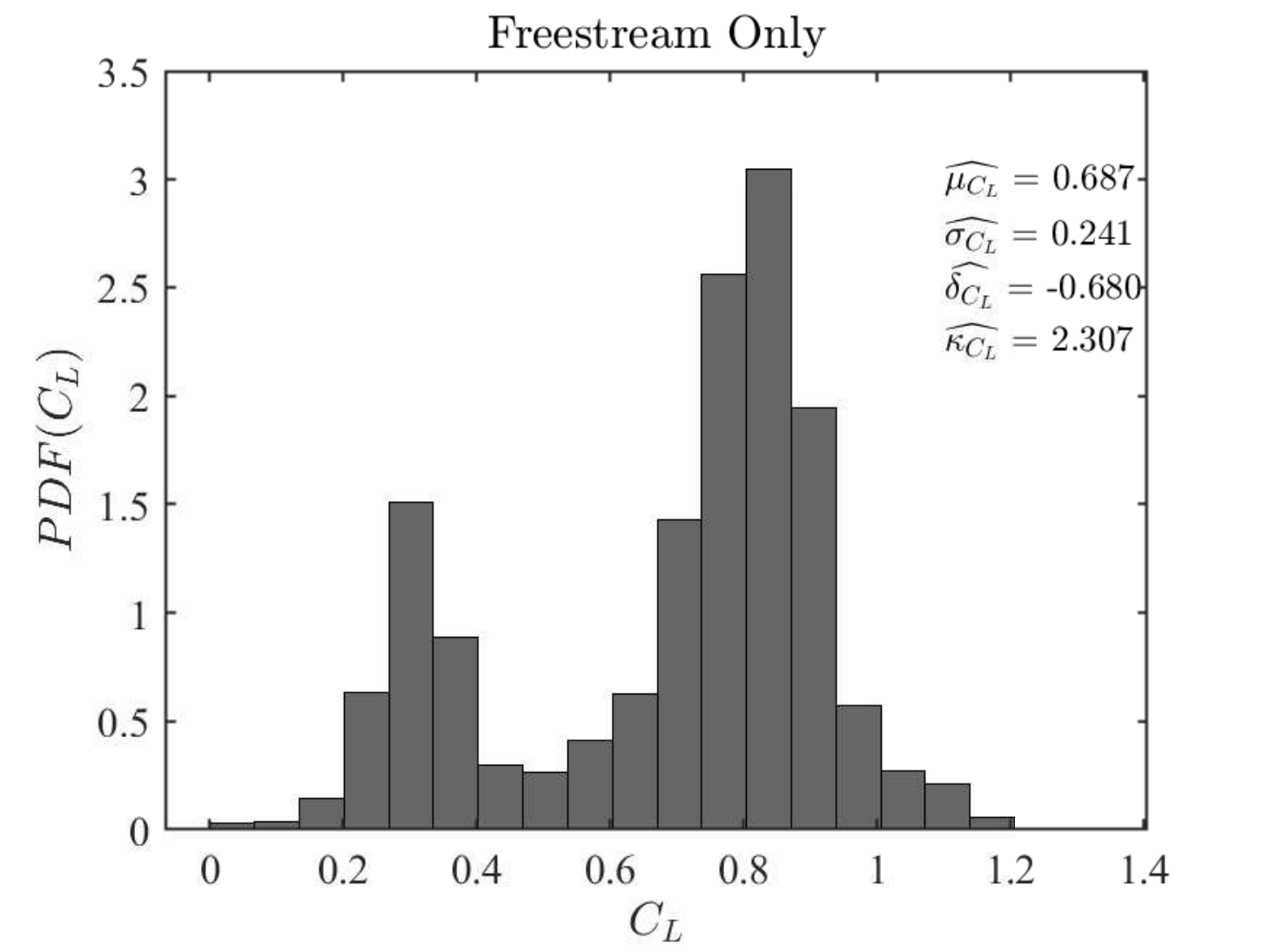}
    \caption{PDF of $C_L$ obtained from $5000$ Monte Carlo simulations for $a_{\infty}  \sim N(4.00 ^{\circ},0.40 ^{\circ})$ and $M_{\infty} \sim N(0.75,0.16)$.} 
    \label{NACA12_FreestreamOnlyPDF} 
\end{figure}

In the first scenario, uncertainty is introduced only in the free-stream conditions, that follow a normal distribution with $a_{\infty}  \sim N(4.0 ^{\circ},0.40 ^{\circ})$ and $M_{\infty} \sim N(0.75,0.16)$. We run  Monte-Carlo simulations with $5000$ samples to quantify the effect of these uncertainties on the QoI. The PDF of $C_L$ is shown in figure \ref{NACA12_FreestreamOnlyPDF}; it has a bimodal distribution with two modes, one with low $C_L$ (in the range $[0.2,0.45]$) and one with high $C_L$ (in the range $[0.6,1.1])$. At the nominal operating conditions, the lift coefficient belongs to the high $C_L$ mode.  The strong deviation from Gaussian distribution makes this a very interesting case.

\begin{figure}[!ht]
\centering
\includegraphics[scale=0.40,clip]{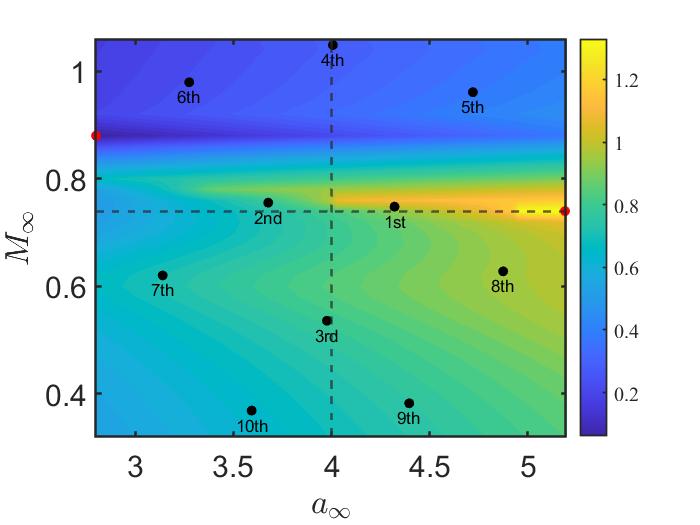}
\caption{Contour plot of $C_L$ in the ($a_{\infty},M_{\infty}$) plane. The QR points for $p = 3$ are also shown and are labeled with their ranking order. The two red dots mark the conditions of minimum and maximum $C_L$ (contours of Mach number are shown in figure \ref{MC_NACA_Freestream_Only_MACH_SHOCKS}). The vertical and horizontal dashed lines denote the nominal free-stream conditions.}
\label{NACA12_FreestreamOnlyContour}
\end{figure}

To shed more light into this behaviour, a contour plot of  $C_L$ in the ($a_{\infty},M_{\infty}$) plane is shown in in figure \ref{NACA12_FreestreamOnlyContour}; the examined rage is $\pm 3 \sigma$ around the nominal values for both stochastic variables. In the plot, we also superimpose the QR sampling points for $p=3$, with the labels indicating the ranking order. As expected, the highest ranked points are closer to the mean. However, contrary to figures \ref{fig:grid_WTS_points} and \ref{fig:WTS_TS_Comparison}, the mean value is not one of the QR points. We have empirically observed that for Gaussian uncertainties the QR points include the mean value only for even $p$, while for odd $p$ the points are placed on either side of the mean. 

As can be seen, the lift coefficient depends strongly on both parameters. It is well known that the variation of $C_L$ with $M_{\infty}$ in the transonic regime is highly non-linear  \cite{Bertin_Cummings_2014}. Visualisation of the flow fields has shown that with increasing the Mach number, while keeping $a_{\infty}$ constant, a shock starts to appear on the suction side of the airfoil at around $M_{\infty} = 0.6$, with area and strength increasing as  $M_{\infty}$ increases. This corresponds to the region of high $C_L$ (green/yellow area on the plot). However, for higher $M_{\infty}$, the lift coefficient  transitions from the high $C_L$ mode to the low $C_L$ mode (blue area). This is associated with a qualitative change in the flow field. 

\begin{figure}[!ht]
     \subfloat[\label{NACA12_FreestreamOnlyContourMACH07}]{\includegraphics[scale=0.16, clip]{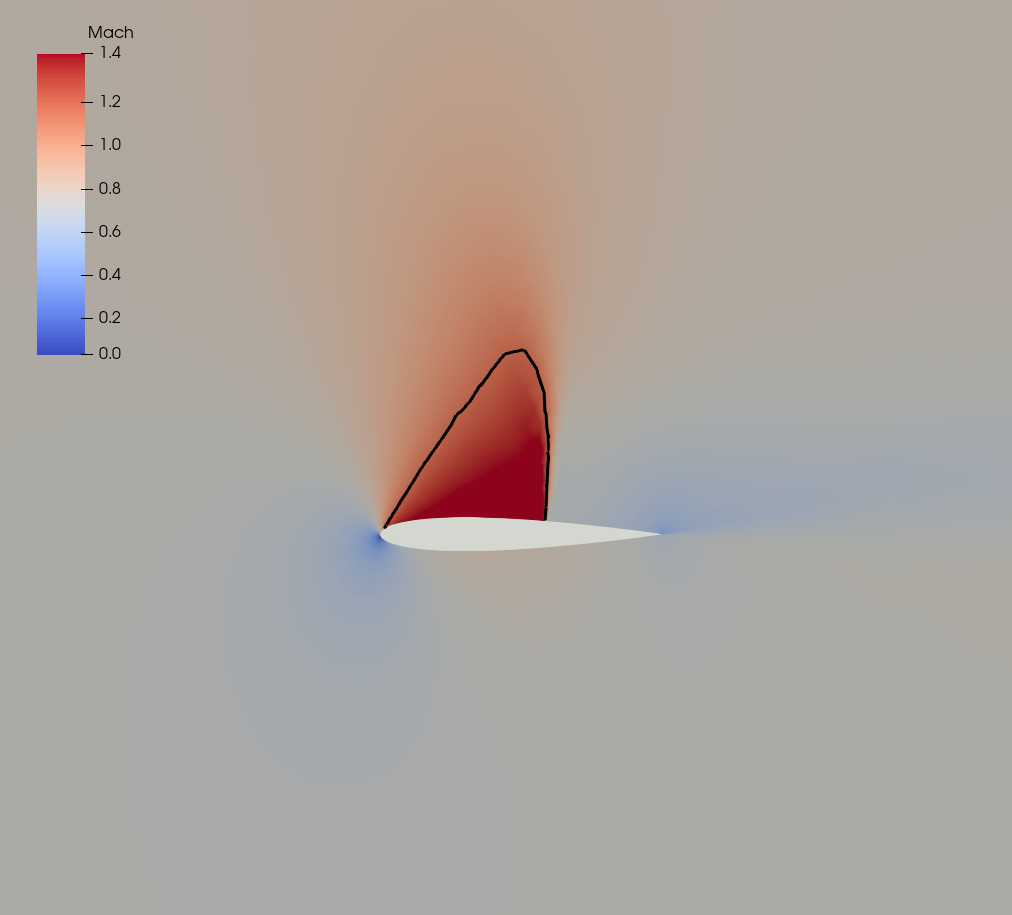}}
     \subfloat[\label{NACA12_FreestreamOnlyContourMACH09}]{\includegraphics[scale=0.16, clip]{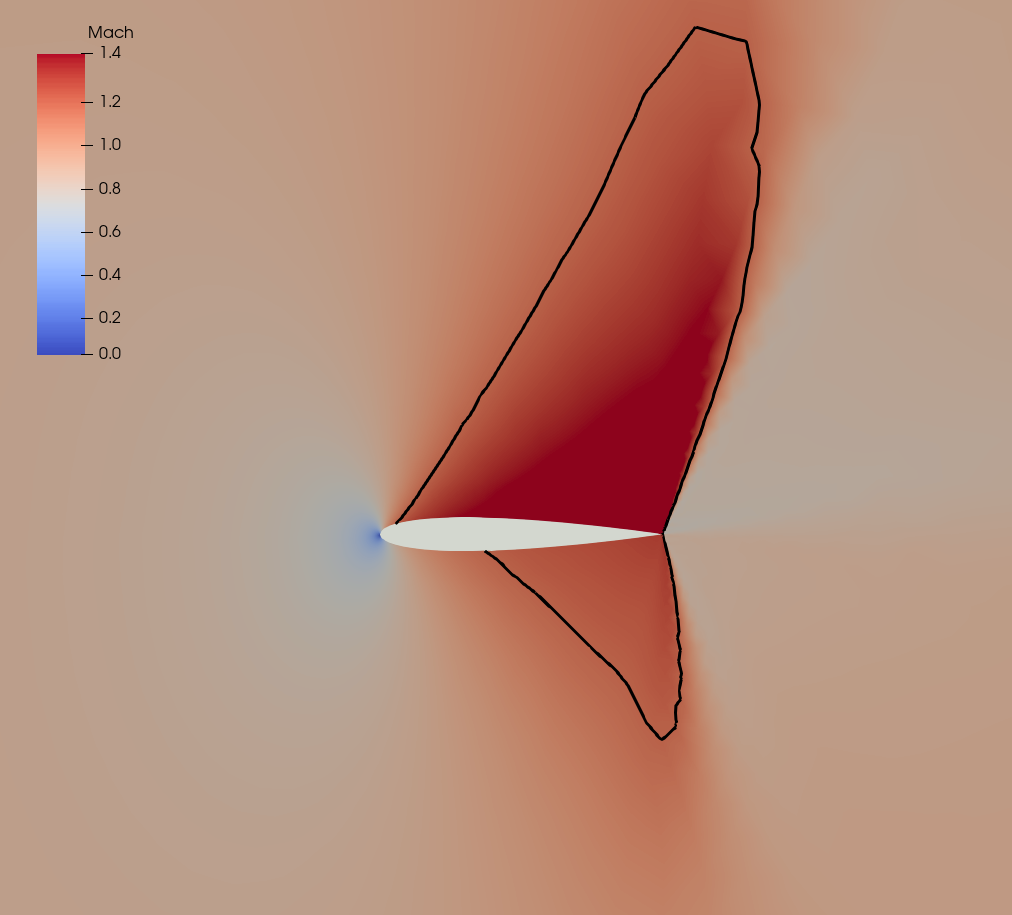}}
	\caption{ Contour plot of $M$ for (a) $M_{\infty } = 0.74$ and $a_{\infty } = 5.2^{\circ}$ and  (b)  $M_{\infty} = 0.88$ and $a_{\infty } = 2.8^{\circ}$. The two conditions are marked with red dots in figure \ref{NACA12_FreestreamOnlyContour}. The black solid line indicates $M = 1$.}
	\label{MC_NACA_Freestream_Only_MACH_SHOCKS}
\end{figure}

In figure \ref{NACA12_FreestreamOnlyContourMACH07} contours of Mach number at $M_{\infty} = 0.74$ and $a_{\infty} = 5.2^{\circ}$ are shown; these conditions (marked with a red dot on the right edge of figure \ref{NACA12_FreestreamOnlyContour}) correspond to the maximum $C_L = 1.35$. It can be seen that only one shock appears on the suction surface of the airfoil. On the other hand, for free-stream conditions $M_{\infty} = 0.88$ and $a_{\infty} = 2.8^{\circ}$ (marked with a red dot on the left edge of figure \ref{NACA12_FreestreamOnlyContour}), a shock appears on both the suction and pressure sides close to the trailing edge, see figure \ref{NACA12_FreestreamOnlyContourMACH09}. In the largest portion of the bottom side the flow is supersonic and the pressure is low, causing the lift to drop to the smallest value of $C_L = 0.06$. More generally, free-stream conditions that lead to the presence of shocks on both sides induce lower $C_L$ coefficients, and are responsible for the low-lift mode of the PDF. On the other hand, conditions that lead to a single shock on the suction side only have a higher lift coefficient and belong to the high-lift mode. This change of the flow pattern explains the bi-modal shape of the PDF. For more details regarding the observed minimum and maximum $C_L$ values and the change of the flow patterns with increasing $M_\infty$ in the transonic flow regime, refer to chapter 9 of \cite{Bertin_Cummings_2014}.  

Comparison of the first 4 moments of the $C_L$ coefficient as computed by Monte Carlo  and se-gPC with $p=3$ are shown in table \ref{NACA12_Comparison_Freestream_Only}. For this case, there are $P+1=10$ unknown coefficients. At each QR point, the value of $C_L$ and the sensitivities $dC_L/dM_{\infty}$, $dC_L/da_{\infty}$ are computed with one direct and one adjoint evaluation, thus we need only 4 QR points (that will provide 12 equations). As can be seen from the table, the mean and standard deviation are computed with less than 1\% error, and even for the kurtosis the error is less than 10\%.  

\begin{table}[!ht]
\centering
\begin{tabular}{|l|l|l|l|l|}
\hline
Method & $\widehat{\mu}_{C_L}$ & $\widehat{\sigma}_{C_L}$ & $\widehat{\delta}_{C_L}$ & $\widehat{\kappa}_{C_L}$ \\ \hline \hline
Monte-Carlo     & $0.687$    & $0.241$   & $-0.680$    & $2.307$    \\ \hline
se-gPC ($p = 3$)  & $0.687$    & $0.240$   & $-0.659$    & $2.511$    \\ \hline
Error ($\%$)  & $<0.01 $    & $0.41$   & $3.08 $    & $8.84$    \\ \hline
\end{tabular}
\caption{First four moments of $C_L$ coefficient as computed by Monte Carlo (with $5000$ samples) and se-gPC  (with $p = 3$).  The Monte Carlo values are considered as reference for the evaluation of percentage error.}
\label{NACA12_Comparison_Freestream_Only}
\end{table}

%%%%%%%%%%%%%%%%%%%%%%%%%%%%%%%%%%%%%%%%%%%%%%%%%%%%%%%%%%
\subsubsection{Geometric and free-stream uncertainties}
%%%%%%%%%%%%%%%%%%%%%%%%%%%%%%%%%%%%%%%%%%%%%%%%%%%%%%%%%%
\label{Transonic Inviscid Naca Freestream and Geometry}

In the second scenario, we consider uncertainties in both the free-stream conditions as well as the geometry. Uncertainties in the geometric shape are present in physical settings and can affect the airfoil performance. To model them, we employ Hicks-Henne bump functions, which have been used to describe shape deformation in aerodynamic optimization and wing design problems \cite{HicksHenne,YANG2018362}.  The parameterized geometry is given by,
\begin{equation}
y(x) = y_{base}(x) + \sum_{i=1}^{m_G} b_i(x),  \quad \quad
b_i(x) = a_i \left [ \sin \left(\pi x^{\frac{log0.5}{logh_i}} \right) \right ]^{t_i},
\label{NACA_Hicks_Henne}
\end{equation}

\noindent where $m_{G}$ is the number of bump functions, $a_i$ is the bump amplitude, $t_i$ controls the width of the bump, and $h_i$ determines the location of the bump maximum.  We consider $m_G=38$ stochastic coefficients $a_i$ that follow normal distribution, $a_i \sim  N(1,0.01)$. The standard deviation is selected to be small; larger values can result in non-physical airfoil shapes. In total, we have $m=m_G+2 = 40$ stochastic parameters, with large variance in the free-stream conditions and small variance for the parameters that describe the geometry. This case mirrors some of the features that are found in robust aerodynamic shape optimization. 
Note that the number of uncertain parameters is not an exaggeration, similar cases in aerodynamic shape optimization use dozens of shape parameters \cite{SHAHROKHI2007443, ShapeOptimizationOnera, ActiveSubspacesOptim}. The open source CFD software $SU2$ \cite{SU2} was used again to simulate the flow and calculate the sensitivities of $C_L$ with respect to the geometric parameters using the adjoint approach. 

\begin{figure}[!ht]
\centering
\includegraphics[scale=0.6, clip]{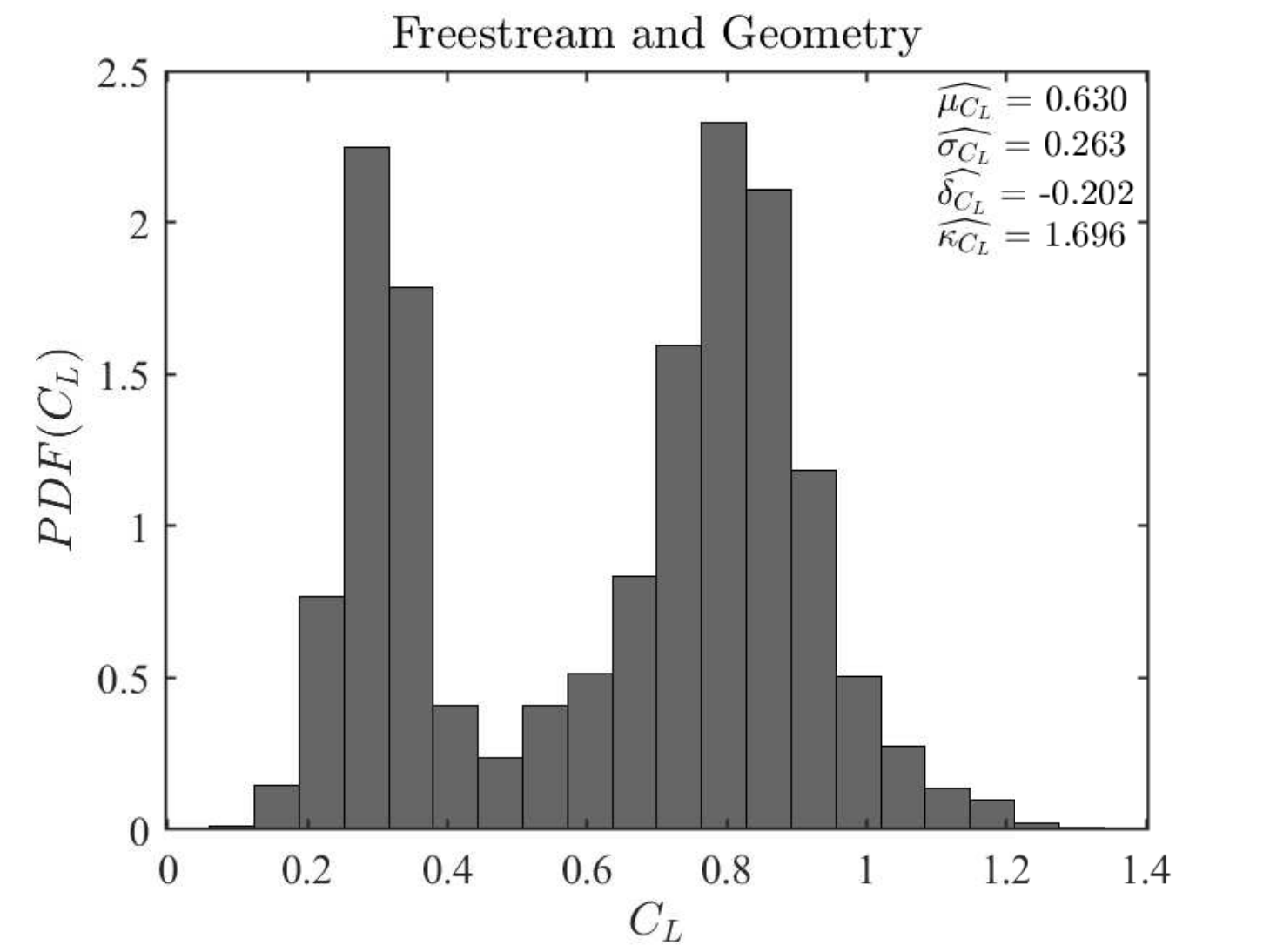}
\caption{PDF of $C_L$ for the NACA12 airfoil produced by $5000$ Monte-Carlo simulations. The $m = 40$ stochastic inputs are $M_{\infty} = N(0.75,0.16)$, $a_{\infty} = N(4^{\circ},0.4^{\circ})$ and $a_i  \sim N(1,0.01)$ with $i = 1,\dots,38$.}
 \label{NACA12_Monte_Carlo}
\end{figure}

The PDF of the lift coefficient $C_L$, obtained from $5000$ Monte Carlo simulations, is shown in figure \ref{NACA12_Monte_Carlo}. The PDF remains bi-modal as in the first scenario, but has  different statistical properties. The performance of the airfoil is significantly affected by the geometric uncertainties. For example, the mean value of the lift is reduced by $8.3\%$ from $0.687$ to $0.630$. The higher order moments are also affected. 

\begin{table}[!ht]
\begin{center}
\begin{tabular}{|l|l|l|l|}
\hline
Method      &  $p = 1$ & $p = 2$    & $p = 3$ \\ \hline \hline
Smolyak    &  
\begin{tabular}{@{}l@{}} $81$ \\ \hline $1.98h$\\ \end{tabular} 
&     
\begin{tabular}{@{}l@{}} $3321$ \\ \hline $80.5h$\\ \end{tabular} &     
\begin{tabular}{@{}l@{}} $91881$ \\ \hline $*$\\ \end{tabular}  
\\ \hline
WLSQ (at QR)    &  
\begin{tabular}{@{}l@{}} $41$ \\ \hline $0.97h$\\ \end{tabular} 
&
\begin{tabular}{@{}l@{}} $861$ \\ \hline $21.2h$\\ \end{tabular}
&
\begin{tabular}{@{}l@{}} $12341$ \\ \hline $*$\\ \end{tabular} 
\\ \hline se-gPC (at QR)    &  
\begin{tabular}{@{}l@{}} $2$ \\ \hline $0.05h$\\ \end{tabular} 
&
\begin{tabular}{@{}l@{}} $42$ \\ \hline $1.04h$\\ \end{tabular} 
&
\begin{tabular}{@{}l@{}} $602$ \\ \hline $15.4h$\\ \end{tabular}     \\ \hline
\end{tabular}
\caption{Number of evaluations (top value in each cell) and CPU hours (bottom value in each cell) required for Smolyak Quadrature, WLSQ and se-gPC at the QR points for different chaos orders $p$. Stars indicate computational cost more expensive than the  Monte Carlo simulations. }
\label{NACA12_Cost_Table}
\end{center}
\end{table}

The computational cost (number of evaluations) of the three methods considered (Smolyak Quadrature, WLSQ, and se-gPC at QR sampling points) is seen in table \ref{NACA12_Cost_Table} for different values of the chaos order $p$. The CPU hours required for each method are also presented (a star indicates that the computational time is longer than the Monte-Carlo simulations). The application of se-gPC results in a very drastic reduction in the number of evaluations and CPU time (one to two orders of magnitude). The difference in performance is much more prominent for this case compared to the previous ones because of the large number of stochastic variables.

 \begin{table}[!ht]
\begin{center}
\begin{tabular}{|l|l|l|l|l|l|l|l|}
\hline
Method  & \begin{tabular}{@{}l@{}} Smolyak \\ $p = 1$\\ \end{tabular} &\begin{tabular}{@{}l@{}} Smolyak \\ $p = 2$\\ \end{tabular} & \begin{tabular}{@{}l@{}} WLSQ \\ $p = 1$\\ \end{tabular} &\begin{tabular}{@{}l@{}} WLSQ \\ $p = 2$\\ \end{tabular}& \begin{tabular}{@{}l@{}} se-gPC \\ $p = 1$\\ \end{tabular} & \begin{tabular}{@{}l@{}} se-gPC \\ $p = 2$\\ \end{tabular} & \begin{tabular}{@{}l@{}} se-gPC \\ $p = 3$\\ \end{tabular}\\ \hline \hline
\begin{tabular}{@{}l@{}} $\widehat{\mu}$ \\ $\epsilon_{\widehat{\mu}}$ \\ \end{tabular}  &     \begin{tabular}{@{}l@{}} $0.638$ \\ \hline $1.26 \%$\\ \end{tabular}     &     \begin{tabular}{@{}l@{}} $0.633$ \\ \hline $0.48 \% $\\ \end{tabular}    &  \begin{tabular}{@{}l@{}} $0.626$ \\ \hline $0.64 \%$\\ \end{tabular}       &   \begin{tabular}{@{}l@{}} $0.627$ \\ \hline $0.48 \%$\\ \end{tabular} & \begin{tabular}{@{}l@{}} $0.635$ \\ \hline $0.79 \%$\\ \end{tabular} & \begin{tabular}{@{}l@{}} $0.633$ \\ \hline $0.48 \%$\\ \end{tabular} &  \begin{tabular}{@{}l@{}} $0.631$ \\ \hline $0.15 \%$\\ \end{tabular}    \\ \hline
\begin{tabular}{@{}l@{}} $\widehat{\sigma}$ \\ $\epsilon_{\widehat{\sigma}}$ \\ \end{tabular}  &       \begin{tabular}{@{}l@{}} $0.259$ \\ \hline $1.52 \%$\\ \end{tabular}     &     \begin{tabular}{@{}l@{}} $0.261$ \\ \hline $0.76 \%$\\ \end{tabular}    &  \begin{tabular}{@{}l@{}} $0.266$ \\ \hline $1.14 \%$\\ \end{tabular}       &   \begin{tabular}{@{}l@{}} $0.264$ \\ \hline $0.38 \%$\\ \end{tabular} & \begin{tabular}{@{}l@{}} $0.266$ \\ \hline $1.14 \%$\\ \end{tabular} & \begin{tabular}{@{}l@{}} $0.265$ \\ \hline $0.76 \%$\\ \end{tabular} &  \begin{tabular}{@{}l@{}} $0.263$ \\ \hline $0.08 \%$\\ \end{tabular}    \\ \hline
\begin{tabular}{@{}l@{}} $\widehat{\delta}$ \\ $\epsilon_{\widehat{\delta}}$ \\ \end{tabular}  &        \begin{tabular}{@{}l@{}} $-$ \\ \hline $-$\\ \end{tabular}     &     \begin{tabular}{@{}l@{}} $-0.182$ \\ \hline $9.9 \%$\\ \end{tabular}    &  \begin{tabular}{@{}l@{}} $-$ \\ \hline $-$\\ \end{tabular}       &   \begin{tabular}{@{}l@{}} $-0.234$ \\ \hline $15.84 \%$\\ \end{tabular} & \begin{tabular}{@{}l@{}} $-$ \\ \hline $-$\\ \end{tabular} & \begin{tabular}{@{}l@{}} $-0.219$ \\ \hline $8.41 \%$\\ \end{tabular} &  \begin{tabular}{@{}l@{}} $-0.209$ \\ \hline $3.46 \%$\\ \end{tabular}    \\ \hline
\begin{tabular}{@{}l@{}} $\widehat{\kappa}$ \\ $\epsilon_{\widehat{\kappa}}$ \\ \end{tabular}  &       \begin{tabular}{@{}l@{}} $-$ \\ \hline $-$\\ \end{tabular}     &     \begin{tabular}{@{}l@{}} $-$ \\ \hline $-$\\ \end{tabular}    &  \begin{tabular}{@{}l@{}} $-$ \\ \hline $-$\\ \end{tabular}       &   \begin{tabular}{@{}l@{}} $-$ \\ \hline $-$\\ \end{tabular} & \begin{tabular}{@{}l@{}} $-$ \\ \hline $-$\\ \end{tabular} & \begin{tabular}{@{}l@{}} $-$ \\ \hline $-$\\ \end{tabular} &  \begin{tabular}{@{}l@{}} $1.799$ \\ \hline $6.07 \%$\\ \end{tabular}    \\ \hline
\end{tabular}
\caption{Comparison of the first four moments of $C_L$ for uncertain free-stream and geometry conditions, as computed by Smolyak Quadrature, WLSQ and se-gPC at the WTS points for different values of the chaos order $p$. Errors are computed with respect to MC with $5000$ sample (reference values are, $\widehat{\mu} = 0.630$, $\widehat{\sigma} = 0.263$, $\widehat{\delta} = -0.202$ and $\widehat{\kappa} = 1.696$. A dash denotes that the moment cannot be computed.}
\label{NACA12_Comp_Table}
\end{center}
\end{table}

Comparison with the Monte Carlo results for the four first moments is shown in table \ref{NACA12_Comp_Table}, with the top value in each cell representing the prediction and the bottom value the $\%$ error with respect to MC with $N = 5000$ evaluations. Se-gPC provides similar accuracy compared to the WLSQ and Smolyak Quadrature, but at a much reduced cost. For example, for $p = 1$, se-gPC provides a very good approximation to the first two statistical moments that is $\sim 20$ times less expensive than WLSQ and $\sim 40$ times less expensive than Smolyak (see table \ref{NACA12_Cost_Table}). As mentioned earlier, for $p = 1$ the cost of the se-gPC is independent of the number of uncertain variables. In the context of aerodynamic design the number of control variables is usually large, and it is useful to be able to get fast estimations of the mean and standard deviation. Larger values of $p$ allow accurate approximations of high order moments also at a much reduced computational cost. For example, the se-gPC produced the most accurate prediction of the skewness and it was the only method that was able to produce an estimation of the kurtosis. 

%%%%%%%%%%%%%%%%%%%%%%%%%%%%%%%%%%%%%%%%%%%%%%%%%%%%%%%%%%%%%%%%%%%%%%%
%%%%%%%%%%%%%%%%%%%%%%%%%%%%%%%%%%%%%%%%%%%%%%%%%%%%%%%%%%%%%%%%%%%%%

\section{Conclusions}
\label{sec:conclusions chapter}
%%%%%%%%%%%%%%%%%%%%%%%%%%%%%%%%%%%%%%%%%%%%%%%%%%%%%%%%%%%%%%%%%%%%%
We propose a method to enrich the Least Squares PCE approach for uncertainty quantification with sensitivity information.  The sensitivities of the quantity of interest with respect to the stochastic variables are computed very efficiently (at the cost of a single evaluation) with the adjoint formulation. Thus at each sampling point, we can obtain a number of equations that is equal to the number of stochastic variables, $m$, with minimal cost. The sampling points to apply the direct and adjoint computations are selected from a large pool of candidate points so as to maximise the determinant of the least squares system matrix. This was achieved by applying pivoted QR decomposition to the measurement matrix. It was shown that these points produce a system with small condition number. The proposed method is called sensitivity-enhanced generalised polynomial chaos, or se-gPC.

For the lowest chaos order $p=1$, se-gPC applied at the minimum number of QR points requires only two evaluations, i.e.\ the cost is independent of the dimension of the stochastic space, $m$. For higher values of $p$, the number of the required evaluations is smaller compared to other methods, such as Smolyak quadrature and the standard LSQ, by a factor equal to $m$. Thus, the savings increase as the dimension of the stochastic space grows. The accuracy of se-gPC was assessed against the standard weighted LSQ as well as Smolyak quadrature and the results were found to be of comparable accuracy for the same chaos order, but were obtained at a much reduced computational cost, between one to two orders of magnitude smaller for the cases examined.  

The se-gPC aims to speed up non-intrusive UQ approaches when adjoint sensitivity information is available. The latter is also useful for optimisation purposes and many open source and commercial packages provide solvers for the adjoint equations. This functionality can be directly exploited for efficient UQ calculations for complex engineering problems with many stochastic variables.    

%%%%%%%%%%%%%%%%%%%%%%%%%%%%%%%%%%%%%%%%%%%%%%%%%%%%%%%%%%%%%%%%%%%%%

\section*{Acknowledgements}
The first author wishes to acknowledge the financial support of the President's Scholarship Award from Imperial College London. 
\bibliography{new}

\end{document}